\documentclass[ijoc,nonblindrev]{informs3_arxiv}

\OneAndAHalfSpacedXII 

\usepackage{natbib}
 \bibpunct[, ]{(}{)}{,}{a}{}{,}%
 \def\bibfont{\small}%

\TheoremsNumberedThrough     
\EquationsNumberedThrough    

\usepackage{mathtools}
\usepackage{graphicx}
\usepackage[caption=false]{subfig}
\usepackage{booktabs} 
\usepackage[dvipsnames]{xcolor}
\usepackage{moresize}
\usepackage{algorithm}
\usepackage{algpseudocode}
\usepackage{bm}

\usepackage{longtable}

\MakeRobust{\Call}

\usepackage[normalem]{ulem}

\usepackage{hyperref}

\usepackage[capitalise,noabbrev]{cleveref}
\usepackage{breqn}
\DeclareMathOperator{\gbtFunc}{GBT}

\DeclareMathOperator{\gbtLeft}{Left}
\DeclareMathOperator{\gbtRank}{rank}
\DeclareMathOperator{\gbtRight}{Right}

\DeclareMathOperator*{\gbtCvx}{cvx}

\DeclareMathOperator*{\weight}{weight}

\DeclareMathOperator*{\cover}{cover}

\DeclareMathOperator*{\diag}{diag}

\newcommand{\gbtVec}[1]{\bm{#1}}
\newcommand{\gbtMat}[1]{\bm{#1}}

\def\addlegendimage{\csname pgfplots@addlegendimage\endcsname}

\Crefname{claim}{Claim}{Claims}
\Crefname{example}{Example}{Examples}
\crefname{appsec}{Appendix}{Appendices}

\colorlet{NotCerulean}{Cerulean!70!Blue}

\begin{document}

\TITLE{Mixed-Integer Convex Nonlinear Optimization with Gradient-Boosted Trees Embedded}
\RUNTITLE{Mixed-Integer Convex Nonlinear Optimization with Gradient-Boosted Trees Embedded}

\ARTICLEAUTHORS{%
	\AUTHOR{Miten Mistry}
	\AUTHOR{Dimitrios Letsios}
	\AFF{Imperial College London, South Kensington, SW7 2AZ, UK.}
	\AUTHOR{Gerhard Krennrich}
	\AUTHOR{Robert M. Lee}
	\AFF{BASF SE, Ludwigshafen am Rhein, Germany.}
	\AUTHOR{Ruth Misener}
	\AFF{Imperial College London, South Kensington, SW7 2AZ, UK.}
}

\RUNAUTHOR{Mistry et al.}


\ABSTRACT{Decision trees usefully represent sparse, high dimensional and noisy data.
Having learned a function from this data, we may want to thereafter integrate the function into a larger decision-making problem,
e.g., for picking the best chemical process catalyst.
We study a large-scale, industrially-relevant mixed-integer nonlinear nonconvex
	optimization problem involving both gradient-boosted trees 
and penalty functions mitigating risk.
This mixed-integer optimization problem with convex penalty terms broadly
	applies to optimizing pre-trained regression tree models.
Decision makers may wish to optimize discrete models to repurpose legacy
	predictive models,
or they may wish to optimize a
	discrete model that accurately represents a data set.
We develop several heuristic methods to find feasible solutions, and an exact,
	branch-and-bound algorithm leveraging structural properties of the
	gradient-boosted trees and penalty functions.
We computationally test our methods on concrete mixture design instance and a chemical catalysis industrial instance. 
}
\maketitle

\section{Introduction}
\label{sec:introduction}

Consider integrating an unknown function into an optimization problem, i.e., without a closed-form formula, but with a data set representing evaluations over a box-constrained feasibility domain.
Optimization in the machine learning literature usually refers to the training procedure, e.g.,\ model accuracy maximization \citep{sra2012,NIPS2012_4522}. 
This paper investigates optimization problems after the training procedure, where the trained predictive model is embedded in the optimization problem.
We consider optimization methods for problems with gradient-boosted tree (GBT) models embedded \citep{10.2307/2699986,elements}.
Advantages of GBTs are myriad 
\citep{Chen:2016:XST:2939672.2939785,NIPS2017_6907}, e.g., they
are robust to scale differences in the training data features, handle both categorical and numerical variables, and can minimize arbitrary, differentiable loss functions. 

\citet{ijcai2018-772} survey approaches for embedding machine learning models as parts of decision-making problems.
We \emph{encode the machine learning model using the native language} \citep{ijcai2018-772}, i.e., in an optimization modeling framework.
Resulting optimization models may be addressed using local \citep{Nocedal2006} or deterministic global \citep{schweidtmann2018global} methods.
The value of global optimization is known in engineering \citep{BOUKOUVALA2016701}, e.g.,\ local minima can lead to infeasible parameter estimation \citep{doi:10.1021/jp0548873} or misinterpreted data \citep{BOLLAS20091768}.
For applications where global optimization is less relevant, we still wish to develop optimization methods for discrete and non-smooth machine learning models, e.g.,\ regression trees.
Discrete optimization methods allow repurposing a legacy model, originally built for prediction, into an optimization framework.
In closely related work, \citet{NIPS2017_7132} train machine learning models to capture the task for which they will be used. This work focusses on generating optimal decisions after the machine learning model is trained.

Our optimization problem incorporates an additional, convex penalty term in the objective. This penalty accounts for risky predicted values where the machine learning model is not well trained due to missing data.
But penalizing distance from the candidate solution to the existing data is not the only reason to add a convex penalty function, e.g., our numerical tests consider an instance with an additional soft constraint.
	\citet{doi:10.1002/aic.690320408} document convex terms common in process engineering: any of those convex nonlinear equations could be incorporated into this framework.
	Another possible application area is in portfolio optimization, e.g., extending the (convex) Markowitz model with cardinality constraint and buy-in threshold constraints \citep{Bienstock1996}.
	Several authors have considered more elaborate extensions, e.g., by integrating uncertainty in the expected return estimate \citep{doi:10.1287/opre.1080.0599} or considering concave transaction costs \citep{Konno2001}.
	But the framework presented in this paper could use GBT models to develop data-driven uncertainty or cost models.

This paper considers a mixed-integer nonlinear optimization problem with convex nonlinearities (convex MINLP). The objective sums a discrete GBT-trained function and a continuous convex penalty function.
We design exact methods computing either globally optimal solutions, or solutions within a quantified distance from the global optimum.
The convex MINLP formulation enables us to solve industrial instances with commercial solvers.
We develop a new branch-and-bound method exploiting both the GBTs combinatorial structure and the penalty function convexity. 
Numerical results substantiate our approach.
The manuscript primarily discusses GBTs, but both the mixed-integer linear programming (MILP) formulation and most of the branch-and-bound methods leverage tree-ensemble structure and can be applied to other tree-ensemble models, e.g., random forests and extremely randomized trees \citep{Breiman2001,Geurts2006}.

	This paper studies a problem that is closely related to \citet{2017arXiv170510883M}.
	Our work differs in that (i) \citet{2017arXiv170510883M} studies a more basic version of our problem formulation (optimizing an objective function derived from tree ensembles whereas our objective also includes a convex penalty) and (ii) our contribution is a specialized branch-and-bound algorithm designed to solve our optimization problem at a large-scale.

\paragraph{Paper organization} 
\Cref{sec:problem_definition} introduces the optimization problem and \Cref{sec:modelling} formulates it as a convex MINLP. 
\Cref{sec:branch_and_bound} describes our branch-and-bound method.
\Cref{sec:case_studies} defines the convex penalty term.
\Cref{sec:numerical_analysis} presents numerical results.
Finally, \cref{sec:discussion} discusses further connections to the literature and \cref{sec:conclusion} concludes.

\section{Background}

This section describes gradient-boosted trees (GBTs) \citep{10.2307/2699986,FRIEDMAN2002367}.
In this paper, GBTs are embedded into the \cref{sec:problem_definition} optimization problem.
GBTs are a subclass of boosting methods \citep{FREUND1995256}.
Boosting methods train many weak learners iteratively that collectively produce a 
	strong learner, where a weak learner is at least better than random guessing.
Each boosting iteration trains a new weak learner against the residual of the 
	previously trained learners by minimizing a loss function.
For GBTs, the weak learners are classification and regression trees \citep{cart_book}.

This paper restricts its analysis to regression GBTs without categorical
	input variables.
A trained GBT function is a collection of binary trees and each of these trees provides 
	its own independent contribution when evaluating at $\gbtVec{x}$.
\begin{definition}\label{def:gbt_function}
	A trained GBT function is defined by sets $(\mathcal{T}, \mathcal{V}_t, \mathcal{L}_{t})$ and 
		values $(i(t,s), v(t,s), F_{t,l})$.  
	The set $\mathcal{T}$ indexes the trees. For a given tree $t\in\mathcal{T}$, $\mathcal{V}_t$ and $\mathcal{L}_t$ index the split and leaf nodes, respectively. At split node $t\in\mathcal{T}$, $s\in\mathcal{V}_t$, $i(t,s)$ and $v(t,s)$ return the split variable and value, respectively. At leaf node $t\in\mathcal{T}$, $l\in\mathcal{L}_{t}$, $F_{t,l}$ is its contribution.
\end{definition}

\begin{figure}[t]
	\FIGURE
	{\includegraphics[width=\linewidth]{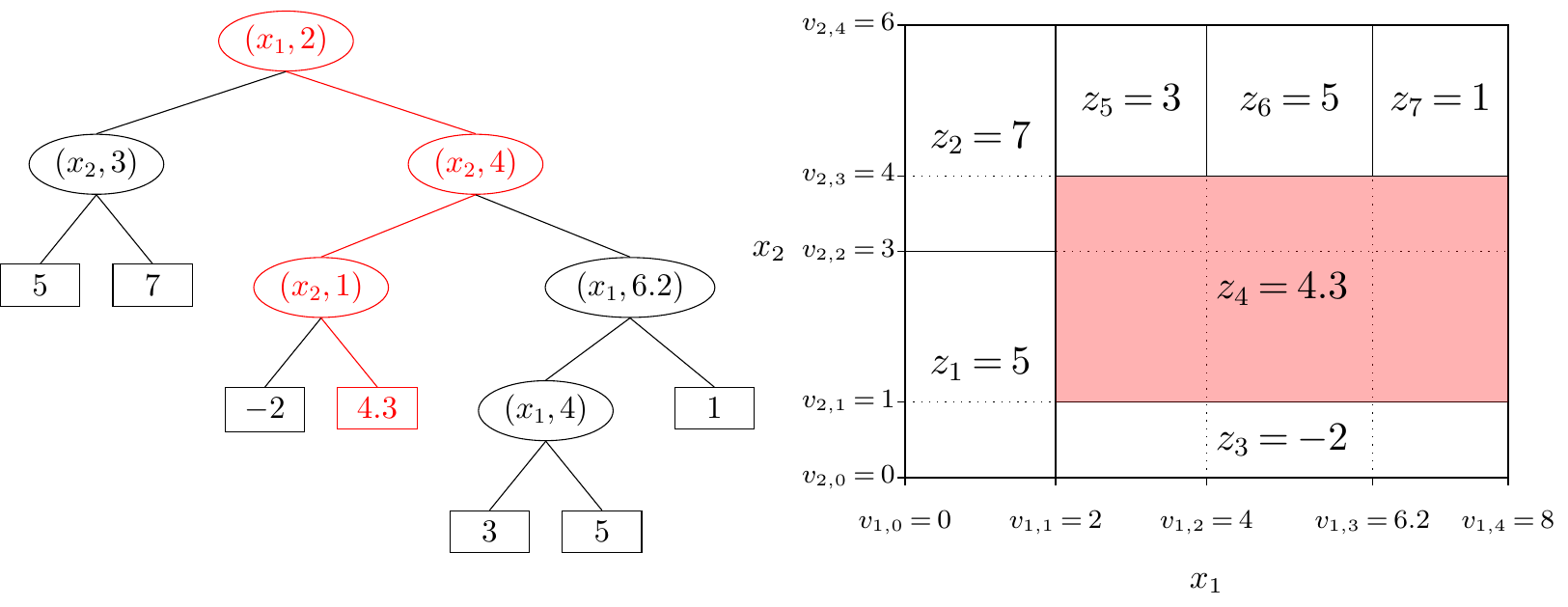}}
	{Gradient boosted tree, see \cref{def:gbt_function}, trained in two dimensions. Left: gradient boosted tree. Right: recursive domain partition defined by tree on left. The highlighted path and region corresponds to the result of evaluating at $\gbtVec{x}=(4.2,2.8)^\top$ as in Example~(1). \label{fig:twod}}
	{}
\end{figure}

\begin{figure}
\FIGURE
	{\includegraphics[width=\linewidth]{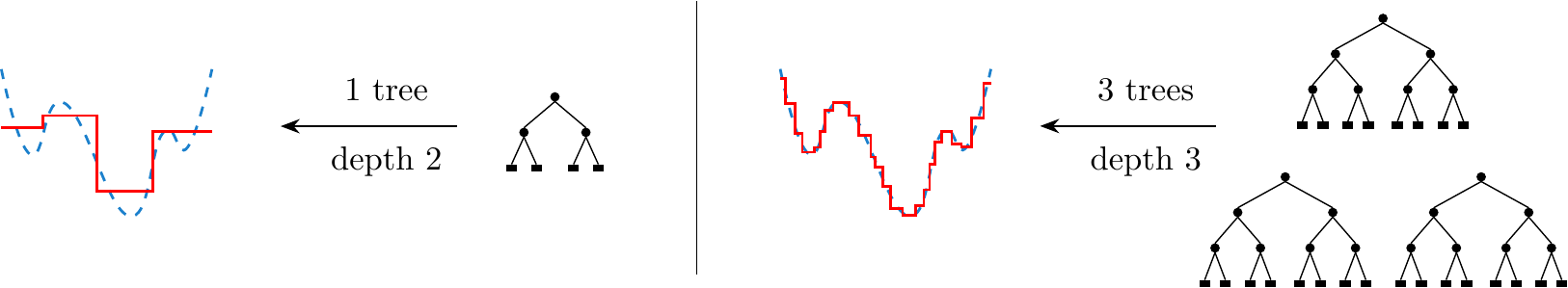}}
	{GBT approximations to the dashed function: 1 tree of depth 2 (left) and 3 trees of depth 3 (right). \label{fig:gbt_example}}
	{}
k\end{figure}

Tree $t\in\mathcal{T}$ evaluates at $\gbtVec{x}$ by following a root-to-leaf path.
Beginning at the root node of $t$, each encountered split node $s\in\mathcal{V}_t$ assesses whether  
	$x_{i(t,s)}<v(t,s)$ or $x_{i(t,s)}\geq v(t,s)$ and follows the left or right child, respectively.
The leaf $l\in\mathcal{L}_{t}$ corresponding to $\gbtVec{x}$ returns $t$'s contribution $F_{t,l}$.
\Cref{fig:twod} shows how a single gradient-boosted tree recursively partitions the domain.
The overall output, illustrated in \Cref{fig:gbt_example}, sums all individual tree evaluations:
\[\gbtFunc(\gbtVec{x})=\sum_{t\in\mathcal{T}}\gbtFunc_t(\gbtVec{x}).\]

\begin{example}\label{ex:gbt}
	Consider a trained GBT that approximates a two-dimensional function with $\mathcal{T}=\{t_1,\dots,t_{|\mathcal{T}|}\}$.
	To evaluate $\gbtFunc(\gbtVec{x})$ where $\gbtVec{x}=(4.2, 2.8)^{\top}$, let $t_1$ be the tree given by \cref{fig:twod}, the highlighted path corresponds to evaluating at $\gbtVec{x}$.
	The root split node query of  $x_1<2$ is false, since  $x_1=4.2$, so we follow the right branch.
	Following this branch encounters another split node. 
	The next query of $x_2<4$ is true, since $x_2=2.8$, so we follow the left branch.
	The final branch reaches a leaf with value 4.3, hence $\gbtFunc_{t_1}(\gbtVec{x})=4.3$. 
	The remaining trees also return a value after making similar queries on $\gbtVec{x}$.
	This results in $\gbtFunc(\gbtVec{x})=\sum_{i=1}^{|\mathcal{T}|}\gbtFunc_{t_i}(\gbtVec{x})=4.3+\sum_{i=2}^{|\mathcal{T}|}\gbtFunc_{t_i}(\gbtVec{x})$.
\end{example}

\section{Optimization Problem}\label{sec:problem_definition}

\begin{table}[t]
	\TABLE
	{Mixed-integer convex programming model sets, parameters and variables.\label{tab:sets_and_params}}
	{\begin{tabular}{cl}
	\toprule
	Symbol & Description\\
	\midrule
	$v^{L}_i$, $v^{U}_i$ & Lower and upper bound of variable $x_i$\\
	$x_i$ & Continuous variable, $i\in\{1,\dots,n\}$
	\\
	$t\in\mathcal{T}$ & Indices of GBTs\\
	$l\in\mathcal{L}_{t}$ & Indices of leaves for tree $t$\\
	$s\in\mathcal{V}_t$ & Indices of split nodes for tree $t$\\
	$m_i$ & Number of variable $x_i$ splitting values \\
	$v_{i,j}$ & Variable $i$'s $j$-th breakpoint, $j\in\{1,\dots,m_i\}$\\
	$F_{t,l}$ & Value of leaf $(t,l)$\\
	$y_{i,j}$ & Binary variable indicating whether variable $x_i<v_{i,j}$\\
	$z_{t,l}$ & Nonnegative variable that activates leaf $(t,l)$\\
	\bottomrule
\end{tabular}
}
	{}
\end{table}

This paper considers box-constrained optimization Problem~(\ref{eq:full_formulation}), an MINLP where the objective sums a convex nonlinear function and a GBT-trained function:
\begin{equation}
	\min_{\gbtVec{v}^L\leq\gbtVec{x}\leq \gbtVec{v}^U}\ 
	\underbrace{
		\gbtCvx(\gbtVec{x})
	}_{\textit{Convex Part}}
	+
	\underbrace{
		\gbtFunc(\gbtVec{x})}_{\textit{GBT Part}}
	,\label[equation]{eq:full_formulation}
\end{equation}
where $\gbtVec{x}=(x_1,\dots,x_n)^{\top}$ is the variable vector.
	$\gbtFunc(\gbtVec{x})$ is the GBT-trained function value at $\gbtVec{x}$.
\Cref{tab:sets_and_params} defines the model sets, parameters and 
		variables.
Problem~(\ref{eq:full_formulation}) is relevant, e.g., when a GBT function has been trained to data but we may trust an optimal solution close to regions with many training points.
A convex penalty term may penalize solutions further from training data.
For instance, consider quality maximization using historical data from a manufacturing process. 
The data may exhibit correlation between two process parameters, e.g., the temperature and the concentration of a chemical additive.
A machine learned model of the system assigns weights to these parameters for future predictions.
Lacking additional information, numerical optimization may produce candidate solutions with temperature and concentration combinations 
that (possibly incorrectly) suggest temperature is responsible for an observed effect.
The convex penalty term 
helps control the optimizer's adventurousness by penalizing deviation from the training data subspace and is parameterized using principal component analysis \citep{Vaswani2018}.
Large values of this risk control term generate conservative solutions.
Smaller penalty values explore regions with greater possible rewards but also additional risk.
Beyond modeling distance to training data, the convex penalty may represent
	additional soft constraints.

A given problem instance may sum independently-trained GBT functions.
Without loss of generality, we equivalently optimize a single GBT function which is the union of all original GBTs.

\section{Mixed-Integer Convex Formulation}
\label{sec:modelling}

Problem~(\ref{eq:full_formulation}) consists of a continuous convex function and a discrete GBT function.
The discrete nature of the GBT function arises from the left/right decisions  
at the split nodes. So we consider a mixed-integer nonlinear program with convex nonlinearities (convex MINLP) formulation.
The main ingredient of the convex MINLP model is a mixed-integer linear programming (MILP)
	formulation of the GBT part which merges with the convex part via a linking constraint.
	The  high level convex MINLP is:
\begin{subequations}
\label{Eq:Optimization_Problem}
	\begin{align}
		\min_{\gbtVec{v}^L\leq\gbtVec{x}\leq\gbtVec{v}^U}\ \ 
		&\gbtCvx(\gbtVec{x}) + [\text{GBT MILP objective}] 
		\label{eq:minlp_ob}\\
		\text{s.t.}\quad &[\text{GBT MILP constraints}],
		\label{eq:minlp_constraints}\\
		&[\text{Variable linking constraints}]\label{eq:minlp_link}.
	\end{align}
	\label{eq:high_level}
\end{subequations}

\subsection{GBT MILP Formulation}

We form the GBT MILP using the \citet{2017arXiv170510883M} approach, which recalls the state-of-the-art in modeling piecewise linear functions \citep{misener-etal:2009,misener-floudas:2010,vielma-etal:2010}.
\citet{VERWER2017368} present an alternative MILP formulation. 
Alternative modeling frameworks include constraint programming \citep{rossi2006handbook,10.1007/978-3-319-18008-3_6} and satisfiability modulo theories \citep{LOMBARDI2017343,MISTRY201898}.

\Cref{fig:twod} shows how a GBT partitions the domain $[\gbtVec{v}^L,\gbtVec{v}^U]$ of $\gbtVec{x}$.
Optimizing a GBT function reduces to optimizing the leaf selection, i.e.,\ 
	finding an optimal interval, opposed to a specific $\gbtVec{x}$ value.
Aggregating over all GBT split nodes produces a vector of ordered breakpoints $v_{i,j}$
	for each $x_i$ variable:
	$v^{L}_i=v_{i,0}<v_{i,1}<\dots<v_{i,m_i}<v_{i,m_i+1}=v^{U}_i$.
Selecting a consecutive pair of breakpoints for each $x_i$ 
defines an interval where the GBT function is constant.
	Each point $x_i\in[v_i^L,v_i^U]$ is either on a breakpoint $v_{i,j}$ or in the interior
	of an interval. 
Binary variable $y_{i,j}$ models whether $x_i< v_{i,j}$ for $i\in[n]=\{1,\dots,n\}$ and $j\in[m_i]=\{1,\dots,m_i\}$.
Binary variable $z_{t,l}$ is 1 if tree $t\in\mathcal{T}$ evaluates at node $l\in\mathcal{L}_t$ and 0 otherwise.
Denote by $\mathcal{V}_t$ the set of split nodes for tree $t$.
Moreover, let $\gbtLeft_{t,s}$ and $\gbtRight_{t,s}$ be the sets of subtree leaf nodes 
rooted in the left and right children of split node $s$, respectively.

MILP Problem~(\ref{eq:ilp_sub}) formulates the GBT \citep{2017arXiv170510883M}.
\Cref{eq:ob_ilp} minimizes the total value of the active leaves.
\Cref{eq:sosz_ilp} selects exactly one leaf per tree.
\Cref{eq:ilp_left_select,eq:ilp_right_select} activates a leaf only if all corresponding splits occur. 
\Cref{eq:ilp_y_order} ensures that if $x_i\leq v_{i,j-1}$, then $ x_i\leq v_{i,j}$.
Without loss of generality, we drop the $z_{t,l}$ integrality constraint because any feasible assignment of $\gbtVec{y}$ specifies one leaf, i.e.,\ a single region in \cref{fig:twod}.
	\begin{subequations}
\begingroup
\allowdisplaybreaks
	\begin{align}
		\min\;
			& \sum_{t\in\mathcal{T}}
						\sum_{l\in\mathcal{L}_{t}}
						F_{t,l}z_{t,l}\label{eq:ob_ilp}
						\\
		\text{s.t.}\;
		&\sum_{l\in\mathcal{L}_{t}}z_{t,l} = 1, &\forall t&\in \mathcal{T}\label{eq:sosz_ilp},\\
		&\sum_{\mathclap{l\in\gbtLeft_{t,s}}}z_{t,l} \leq y_{i(s),j(s)}, &\forall t&\in\mathcal{T}, s\in 
		\mathcal{V}_{t} \label{eq:ilp_left_select},\\
		&\sum_{\mathclap{l\in\gbtRight_{t,s}}}z_{t,l} \leq 1-y_{i(s),j(s)}, &\forall t&\in\mathcal{T}, s\in \mathcal{V}_{t}\label{eq:ilp_right_select},\\
		&y_{i,j}\leq y_{i,j+1}, &\forall i&\in[n],\,j\in[m_i-1]\label{eq:ilp_y_order},\\
			&y_{i,j}\in\{0,1\}, &\forall i&\in[n],\, j\in [m_i], \\
			&z_{t,l}\geq0, &\forall t&\in\mathcal{T},\,l\in\mathcal{L}_{t}.
	\end{align}
\endgroup
		\label{eq:ilp_sub}%
	\end{subequations}%
\subsection{Linking Constraints}

\Cref{eq:link1,eq:link2} relate the continuous $x_i$ variables, from the original Problem~(\ref{eq:full_formulation}) definition, to the binary $y_{i,j}$ variables:
\begin{subequations}
	\begin{align}
		x_i&\geq v_{i,0} + \sum_{j=1}^{m_i}(v_{i,j} - v_{i,j-1})(1-y_{i,j}),\label{eq:link1}\\
		x_i&\leq v_{i,m_i+1}+ \sum_{j=1}^{m_i}(v_{i,j} - v_{i,j+1})y_{i,j},\label{eq:link2}
	\end{align}%
	\label{eq:link}%
\end{subequations}%
	for all $i\in[n]$.
We express the linking constraints using non-strict inequalities to
	avoid computational issues when optimizing with strict inequalities.
Combining \cref{eq:high_level,eq:ilp_sub,eq:link} defines the mixed-integer nonlinear program with convex nonlinearities (convex MINLP) formulation to Problem~(\ref{eq:full_formulation}). \cref{app:full_micp} lists the complete formulation.

\subsection{Worst Case Analysis}\label{subsec:analysis}
The difficulty of Problem~(\ref{eq:full_formulation}) is primarily 
justified by the fact that optimizing a GBT-trained function, i.e.,\ Problem~(\ref{eq:ilp_sub}), 
is an NP-hard problem \citep{2017arXiv170510883M}.
This section shows that the number of continuous
variable splits and tree depth affects complete enumeration. 
These parameters motivate the branching scheme in our branch-and-bound algorithm.

In a GBT ensemble, each continuous variable $x_i$ is associated with $m_i+1$ intervals (splits).
Picking one interval $j\in\{1,\dots,m_i+1\}$ for each $x_i$ sums to a total of $\prod_{i=1}^n(m_i+1)$ distinct combinations.
A GBT-trained function evaluation selects a leaf from each tree.
But not all leaf combinations are valid evaluations.
In a feasible leaf combination where one leaf enforces $x_i<v_1$ and another enforces $x_i\geq v_2$, it must be that $v_2<v_1$.
Let $d$ be the maximum tree depth in $\mathcal{T}$.
Then the number of leaf combinations is upper bounded by $2^{d|\mathcal{T}|}$.  
Since the number of feasibility checks for a single combination is 
	$\tfrac{1}{2}|\mathcal{T}|(|\mathcal{T}|-1)$,
an upper bound on the total number of feasibility checks is
	$2^{d|\mathcal{T}|-1}|\mathcal{T}|(|\mathcal{T}|-1)$.
So the worst case performance
of an exact method improves as the number of trees decreases.

\section{Branch-and-Bound Algorithm}
\label{sec:branch_and_bound}

This section designs an exact branch-and-bound (B\&B) approach.
Using a divide-and-conquer principle, B\&B forms a tree of subproblems and 
searches the domain of feasible solutions.
Key aspects of B\&B are: (i) rigorous lower (upper) bounding methods 
	for minimization (maximization) subproblems, (ii) branch variable and value selection, and (iii) 
	feasible solution generation.
In the worst case, B\&B enumerates all solutions, but generally it avoids 
	complete enumeration by pruning subproblems, i.e.,\ removing infeasible subproblems 
	or nodes with lower bound exceeding the best found feasible solution \citep{MORRISON201679}.
This section exploits 
\emph{spatial branching} that splits on continuous variables 
\citep{belotti_kirches_leyffer_linderoth_luedtke_mahajan_2013}.
\Cref{Table:Nomenclature} in \cref{Appendix:Nomenclature} defines the symbols in this section.

\subsection{Overview}

\begin{algorithm}[t]
\caption{Branch-and-Bound (B\&B) Algorithm Overview}
\label{alg:bb_overview}
\begin{algorithmic}[1]
	\State $S=[\bm{L},\bm{U}]\leftarrow[\bm{v}^L,\bm{v}^U]$
	\State $b^{\text{cvx},S}\leftarrow\Call{ConvexBound}{S}$
	\Comment{\cref{Lem:Global_Lower_Bound}, \cref{Section:Global_Lower_Bound}}
	\State $P_{\text{root}}\leftarrow\Call{RootNodePartition}{N}$ 
	\Comment{\cref{subsubsec:gbt_bound}}
	\State $b^{\text{GBT},S,P_{\text{root}}}\leftarrow\Call{GbtBound}{S,P_{\text{root}}}$
	\Comment{\cref{lem:gbt_lower_bound}, \cref{subsubsec:gbt_bound}}
	\State $B\leftarrow\Call{BranchOrdering}$ 
	\Comment{\cref{subsec:preprocess}}
	\State $Q=\{S\}$
	\While{$Q\neq\emptyset$}
		\State Select $S\in Q$
		\If {$S$ is not leaf}
		\State $S'\gets S$
		\Repeat
			\State $S',(x_i,v)\leftarrow\Call{StrongBranch}{S',B}$
			\Comment{\cref{alg:strong_branch}, \cref{subsec:strong_branch}}
		\Until{strong branch not found}
		\If{$S'$ is not leaf}
		\State$(S_{\text{left}},S_{\text{right}})\leftarrow\Call{Branch}{S',(x_i,v)}$
		\State $P$: tree ensemble partition of node $S$
		\State $P'\leftarrow\Call{PartitionRefinement}{P}$ 
		\Comment{\cref{Alg:Partition_Refinement}, \cref{subsubsec:gbt_bound}}
		\State $b^{\text{GBT},S',P'}\leftarrow\Call{GbtBound}{S',P'}$
		\Comment{\cref{lem:gbt_lower_bound}, \cref{subsubsec:gbt_bound}}		
		\For {$S_{\text{child}}\in\{S_{\text{left}},S_{\text{right}}\}$}
		\If {$S_{\text{child}}$ cannot be pruned}
		\Comment{\cref{Section:Node_Pruning}}
		\State $Q\leftarrow Q\cup\{S_{\text{child}}\}$
		\EndIf
		\EndFor
		\EndIf
		\EndIf
		\State $Q\leftarrow Q\setminus\{S\}$
	\EndWhile
\end{algorithmic}
\end{algorithm}

B\&B \cref{alg:bb_overview} spatially branches over the $[\gbtVec{v}^L,\gbtVec{v}^U]$ domain.
It selects a variable $x_i$, a point $v$ and splits interval $[v_i^L,v_i^U]$ into intervals 
$[v_i^L,v]$ and $[v,v_i^U]$.
Each interval corresponds to an independent subproblem and a new B\&B node.
To avoid redundant branches, all GBT splits define the B\&B branching points.
At a given node, denote the reduced node domain 
by $S=[\gbtVec{L},\gbtVec{U}]$.
\Cref{alg:bb_overview} solves Problem~(\ref{eq:full_formulation}) by relaxing the \cref{eq:link} linking constraints
 and thereby separating the convex and GBT parts.
Using this separation, \cref{alg:bb_overview} computes corresponding bounds 
$b^{\text{cvx},S}$ and $b^{\text{GBT},S,P}$ independently, 
where the latter bound requires a tree ensemble partition $P$ initialized at the root node
and dynamically refined at each non-root node.

\Cref{alg:bb_overview} begins by constructing the root node, 
computing a global lower bound,
and determining a global ordering of all branches (lines 1--5). 
A given iteration: (i) extracts a node $S$ from the unexplored node set $Q$, 
(ii) strong branches at $S$ to cheaply identify branches that tighten the domain resulting in node $S'$, 
(iii) updates the GBT lower bound at $S'$,
(iv) branches to obtain the child nodes $S_{\text{left}}$ and $S_{\text{right}}$, 
(v) assesses if each child node $S_{\text{child}}\in\{S_{\text{left}},S_{\text{right}}\}$ 
may now be pruned and, if not,
(vi) adds $S_{\text{child}}$ to the unexplored node set $Q$ (lines 8--25).

The remainder of this section is structured as follows.
\Cref{subsec:bb_lb} lower bounds Problem~(\ref{eq:full_formulation}).
\cref{Section:Branching} introduces a GBT branch ordering 
and leverages strong branching for cheap node pruning.
\cref{subsec:feasibility} discusses heuristics for computing efficient upper bounds.

\subsection{Lower Bounding}\label{subsec:bb_lb}
\subsubsection{Global lower bound}
\label{Section:Global_Lower_Bound}

The convex MINLP Problem~(\ref{Eq:Optimization_Problem}) objective function consists of a \emph{convex} (penalty) part and a \emph{mixed-integer linear} (GBT) part. 
Lemma \ref{Lem:Global_Lower_Bound} computes a lower bound on the problem by handling the convex and GBT parts independently. 

\begin{lemma}
\label{Lem:Global_Lower_Bound}
Let $S=[\gbtVec{L},\gbtVec{U}]\subseteq[\gbtVec{v}^L,\gbtVec{v}^U]$ be a sub-domain of optimization Problem~(\ref{Eq:Optimization_Problem}). 
Denote by $R^S$ the optimal objective value, i.e., the tightest relaxation, over the sub-domain $S$. 
Then, it holds that $R^S\geq\hat{R}^S$, where:
\begin{equation*}
\hat{R}^S = 
\underbrace{\left[\min_{\gbtVec{x}\in S}\ \gbtCvx(\gbtVec{x}) \right]}_{b^{\gbtCvx, S}}\\
+\underbrace{\left[\min_{\gbtVec{x}\in S}\sum_{t\in\mathcal{T}}\gbtFunc_t(\gbtVec{x})\right]}_{b^{\gbtFunc,S,*}}.
\end{equation*}
\label{lem:global_lower_bound}
\end{lemma}
\proof{Proof}
Let $\gbtVec{x}^*=\argmin_{\gbtVec{x}\in S}\{\gbtCvx(\gbtVec{x}) + \gbtFunc(\gbtVec{x})\}$ and observe 
that $\gbtCvx(\gbtVec{x}^*)\geq b^{\gbtCvx, S}$ and $\gbtFunc(\gbtVec{x}^*)\geq b^{\gbtFunc,S,*}$.
\Halmos
\endproof

\medskip

We may compute $\hat{R}^S$ by removing the \cref{eq:link} linking constraints and solving the mixed-integer model
consisting of \cref{eq:high_level,eq:ilp_sub}.
Computationally, the \cref{lem:global_lower_bound} separation leverages 
efficient algorithms for the convex part and commercial codes 
for the MILP GBT part.
\Cref{lem:global_lower_bound} treats the two Problem~(\ref{eq:full_formulation}) objective terms independently,  
i.e., $\hat{R}^S$ separates the convex and GBT parts.
The \cref{lem:global_lower_bound} separation, while loose at the root node, may be leveraged to discard regions that are dominated by an objective term. Our approach resembles exact algorithms for multiobjective optimization \citep{fernandez2009,nieblinggow,doi:10.1137/18M1169680}. 
An alternative approach, e.g., in line with augmented Lagrangian methods for stochastic optimization \citep{bertsekas2014constrained}, would not separate the convex penalty term as in \cref{lem:global_lower_bound}, but rather tighten the lower bound by integrating the convex penalty and GBTs.
This would be an interesting alternative, but would eliminate the possibility of the strong branching method used in \cref{subsec:strong_branch}.

\subsubsection{GBT Lower Bound}\label{subsubsec:gbt_bound}
While we may efficiently compute $b^{\gbtCvx, S}$ \citep{Boyd2004}, 
deriving $b^{\gbtFunc, S, *}$ is $\mathcal{NP}$-hard \citep{2017arXiv170510883M}.
With the aim of tractability, we calculate a relaxation of $b^{\gbtFunc, S, *}$.
\Cref{lem:gbt_lower_bound} lower bounds Problem~(\ref{eq:ilp_sub}), i.e., the GBT part of Problem~(\ref{Eq:Optimization_Problem}),  by partitioning the GBT ensemble into a collection of smaller ensembles. 

\begin{lemma}
\label{lem:gbt_lower_bound}
Consider a sub-domain $S=[\gbtVec{L},\gbtVec{U}]\subseteq[\gbtVec{v}^L,\gbtVec{v}^U]$ of the optimization problem.
Let $P = \{\mathcal{T}_1,\ldots,\mathcal{T}_k\}$ be any partition of $\mathcal{T}$, i.e., 
$\cup_{i=1}^k\mathcal{T}_i=\mathcal{T}$ and $\mathcal{T}_i\cap\mathcal{T}_j=\emptyset$ $\forall 1\leq i<j\leq k$. 
Then, it holds that $b^{\gbtFunc,S,*}\geq b^{\gbtFunc,S,P}$, where:
\begin{equation*}
b^{\gbtFunc,S,P} = \sum_{\mathcal{T}'\in P}
\left[\min_{\gbtVec{x}\in S}\left\{\sum_{t\in\mathcal{T}'}\gbtFunc_t(\gbtVec{x})\right\}\right].
\end{equation*}
\end{lemma}
\proof{Proof}
When evaluating $\gbtFunc(\gbtVec{x})$ at a given $\gbtVec{x}$, each tree $t\in\mathcal{T}$ provides its own independent contribution $\gbtFunc_t(\gbtVec{x})$, i.e.,\ a single leaf.
A feasible selection of leaves has to be consistent with respect to the GBT node splits, i.e.,\ if one leaf splits on $x_i<v_1$ and another splits on $x_i\geq v_2$ then $v_1 > v_2$.	
Relaxing this consistency requirement by considering a partition $P$ of $\mathcal{T}$ derives the lower bounds $b^{\text{GBT}, S,P}$ for any partition $P$.
\Halmos
\endproof

\medskip

\paragraph{Root Node Partition}
B\&B \cref{alg:bb_overview} chooses an initial root node partition $P_{\text{root}}$ with subsets of size $N$ and calculates the associated \cref{lem:gbt_lower_bound} lower bound. 
Section~\ref{sec:numerical_analysis} numerically decides the partition size $N$ for the considered instances.
The important factors for a subset size $N$ are the tree depth, the number of continuous 
	variable splits and their relation with the number of binary variables.

\paragraph{Non-Root Node Partition Refinement}
Any non-root B\&B node has reduced domain $\gbtVec{x}\in S=[\gbtVec{L},\gbtVec{U}]\subset[\gbtVec{v}^L,\gbtVec{v}^U]$.
B\&B \cref{alg:bb_overview} only branches on GBT node splits, so modeling the reduced 
	domain $S$ in MILP Problem~(\ref{eq:ilp_sub}) is equivalent to setting $y_{i,j}=0$ or $y_{i,j}=1$ for any $y_{i,j}$ 
	that corresponds to $x_i\leq L_i$ or $x_i\geq U_i$, respectively.
	Reducing the box-constrained domain at the node level equates to reducing the GBT instance size. In particular, we may reduce the number and height of trees by assigning fixed variable values and cancelling redundant constraints \citep{2017arXiv170510883M}.

Assume that, at some non-root node with domain $S$, the algorithm is about to update $b^{\text{GBT}, S'\!, P'}$ which was calculated at the parent node with domain $S'\supset S$. 
Fixing binary variables $y_{i,j}$ subject to domain $S$ reduces the worst case enumeration cost of calculating $b^{\text{GBT}, S,P'}$.
The GBT lower bound may further improve at $S$ by considering an alternative partition $P$ such that $|P|<|P'|$, i.e., reducing the number of subsets.
However, reducing the number of subsets has challenges because: (i) choosing any partition $P$ does not necessarily guarantee $b^{\text{GBT},S,P}\geq b^{\text{GBT},S'\!,P'}$, and (ii) a full \cref{lem:gbt_lower_bound} calculation of $b^{\text{GBT},S,P}$ may still be expensive when considering the cumulative time across all B\&B nodes.
Refinability \cref{def:refinability} addresses the choice of $P$ such that $b^{\text{GBT},S,P}\geq b^{\text{GBT},S'\!,P'}$.

\begin{definition}\label{def:refinability}
Given two partitions $P'$ and $P''$ of set $\mathcal{T}$, we say that $P'$ \emph{refines} $P''$ if and only if $\forall \mathcal{T}'\in P',\; \exists \mathcal{T}''\in P''$ such that $\mathcal{T}'\subseteq \mathcal{T}''$.
This definition of refinement implies a partial ordering between different partitions of $\mathcal{T}$.
We express the refinement relation by $\preceq$, i.e., $P'\preceq P''$ if and only if $P'$ refines $P''$.
\end{definition}
\begin{example}
		Let $P=\{\{1,2,3\},\{4,5\}\}$, $P'=\{\{1\},\{2\},\{3\},\{4\},\{5\}\}$ and $P''=\{\{1,2\}, \{3,4,5\}\}$ be partitions of $\{1,\dots,5\}$.
	Here $P'$ refines $P$ since every subset in $P'$ is a subset of one of the $P$ subsets. Similarly $P'$ refines $P''$.
	Partition $P$ does not refine $P''$ nor does $P''$ refine $P$.
\end{example}

\cref{lem:lower_bound_improve} allows bound tightening by partition refinements. Its proof is similar to \cref{lem:gbt_lower_bound}.

\begin{lemma}
\label{lem:lower_bound_improve}
Let $P$ and $P'$ be two partitions of $\mathcal{T}$. If $P'\preceq P$, then $b^{\gbtFunc,P'}\leq b^{\gbtFunc,P}$.
\end{lemma}

In general, for two partitions $P$ and $P'$, we do not know a priori which partition results 
in a superior GBT lower bound.
However, by \cref{lem:lower_bound_improve}, $P'$ refining $P$ suffices for $b^{\text{GBT},P} \geq b^{\text{GBT},P'}$.
Therefore, given partition $P'$ for the parent node, constructing $P$ for the child node $S$ by unifying subsets of $P'$ will not result in inferior lower bounds.

\Cref{Alg:Partition_Refinement} improves $b^{\text{GBT}, S'\!, P'}$ at node $S$ 
by computing a refined partition $P$.
Suppose that $P'= \{\mathcal{T}_1, \dots, \mathcal{T}_k\}$.
Each GBT ensemble subset $\mathcal{T}'\in P'$ corresponds to a smaller subproblem 
with $n^{\mathcal{T}',S}$ leaves ($z_{t,l}$ variables) over the domain $S$.
Initially, \cref{Alg:Partition_Refinement} sorts the subsets of $P'$ in non-decreasing order of 
$n^{\mathcal{T}',S}$.
Then, it iteratively takes the union of consecutive pairs and calculates the associated 
lower bound, i.e., the first calculation is for $b^{\text{GBT}, S, \{\mathcal{T}_1\cup \mathcal{T}_2\}}$, 
the second is for $b^{\text{GBT}, S,\{\mathcal{T}_3\cup \mathcal{T}_4\}}$ and so forth.
The iterations terminate when all unions have been recalculated, or at user defined time limit $q$
resulting in two sets of bounds: those that are combined and recalculated, and those that remain unchanged.
Assuming that the final subset that is updated has index $2l$, the 
	new partition of the trees at node $S$ is
		$P=\{\mathcal{T}_1\cup \mathcal{T}_2, \dots, \mathcal{T}_{2l-1}\cup \mathcal{T}_{2l}, 
		\mathcal{T}_{2l+1},\dots, \mathcal{T}_k\}$
	with GBT bound
		$b^{\text{GBT},S,P}=\sum_{i=1}^l b^{\text{GBT}, S,\{\mathcal{T}_{2i-1}\cup \mathcal{T}_{2i}\}} + 
		\sum_{i=2l+1}^k b^{\text{GBT},S'\!,\{\mathcal{T}_i\}}$.
The second sum is a result of placing time limit $q$ on updating the GBT lower bound.
Time limit $q$ maintains a balance between searching and bounding.
Unifying any number of subsets satisfies \cref{lem:lower_bound_improve}, but  
	\cref{Alg:Partition_Refinement} unifies pairs to keep the resulting subproblems manageable.
	One may speed up our lower bounding procedure by reducing the height of the GBTs, thus relaxing feasibility, and converting each partition subset $\mathcal{T}_k$ solution into a feasible one for $\mathcal{T}_k$ using the \citet{2017arXiv170510883M} split generating procedure for fixing violated constraints.

\begin{algorithm}[t]
	\caption{Non-Root Node Partition Refinement}
	\label{Alg:Partition_Refinement}
	\begin{algorithmic}[1]
	\State $P'$: parent node partition
	\State Sort $P'=\{\mathcal{T}_1,\ldots,\mathcal{T}_k\}$ so that
	$n^{\mathcal{T}_1} \leq \ldots \leq n^{\mathcal{T}_k}$
	\State $P\gets\emptyset$
	\State $i=1$
	\While {$i<\lfloor n/2\rfloor$ and the time limit is not exceeded}
	\State $P\leftarrow P\cup\{\mathcal{T}_{2i-1}\cup\mathcal{T}_{2i}\}$
	\State $i\leftarrow i+1$
	\EndWhile
	\State $P\leftarrow P\cup \{\mathcal{T}_j\in P':j>i\}$
	\State \Return $P$
\end{algorithmic}
\end{algorithm}

\subsubsection{Node Pruning}
\label{Section:Node_Pruning}

In the B\&B algorithm, each node can access: (i) the current best 
	found feasible objective $f^*$, (ii) a lower bound on the convex penalties $b^{\text{cvx}, S}$, 
	and (iii) a lower bound on the GBT part $b^{\text{GBT}, S}$.
The algorithm prunes node $S$ if:
\begin{equation}
	b^{\text{cvx},S} + b^{\text{GBT},S} > f^*,\label{eq:node_prune}
\end{equation}
i.e., if all feasible solutions in $S$ have objective inferior to $f^*$.

\subsection{Branching}
\label{Section:Branching}

\subsubsection{Branch Ordering}
\label{subsec:preprocess}

Next branch selection is a critical element of B\&B \cref{alg:bb_overview}.
Each branch is a GBT split $(x_i,v)$ choice and eliminates a certain number of GBT leaves.
Branching with respect to a GBT split that covers a larger number of 
	leaves may lead to a smaller number of subsequent B\&B iterations by reducing the GBT size.

	Selecting a $(x_i, v)$ split that most improves the GBT lower bound is challenging as it may require solving multiple expensive MILPs.
	So, we heuristically approximate objective improvement by quantifying splits that (i) occur often among all trees and (ii) influence a larger number of leaves in participating trees.
	Let $r((x_i, v), t)$  and $\cover(s,\,t)$ return the set of nodes in tree $t$ that split on $(x_i, v)$ and the set of leaves that node $s\in \mathcal{V}_t$ covers, respectively. 
	We initialize pseudocosts by weighting the $(x_i, v)$ splits:
	\begin{subequations}
		\begin{align}
			\weight((x_i, v), t) &= |\mathcal{L}_t|^{-1}\sum_{\mathclap{s\in r((x_i, v), t)}}|\cover(s, t)|,\label{eq:cover_weight_tree}\\
			\weight((x_i, v), \mathcal{T}) &= \sum_{t\in\mathcal{T}}\weight((x_i, v), t).\label{eq:cover_weight_sum}
		\end{align}%
		\label{eq:weight_func_cover}%
	\end{subequations}%
	\Cref{eq:cover_weight_tree} weights $(x_i, v)$ as the fraction of leaves covered by nodes splitting on $(x_i, v)$ in tree $t$.
Recall that $|\mathcal{L}_t|$ is the number of leaf nodes in tree $t$.
	\Cref{eq:cover_weight_sum} sums all weights calculated by \cref{eq:cover_weight_tree} for split $(x_i,v)$ in each tree $t\in\mathcal{T}$.
	The splits are sorted in non-increasing order their pseudocosts.
\begin{figure}[t]
	\FIGURE
	{\includegraphics[width=0.8\linewidth]{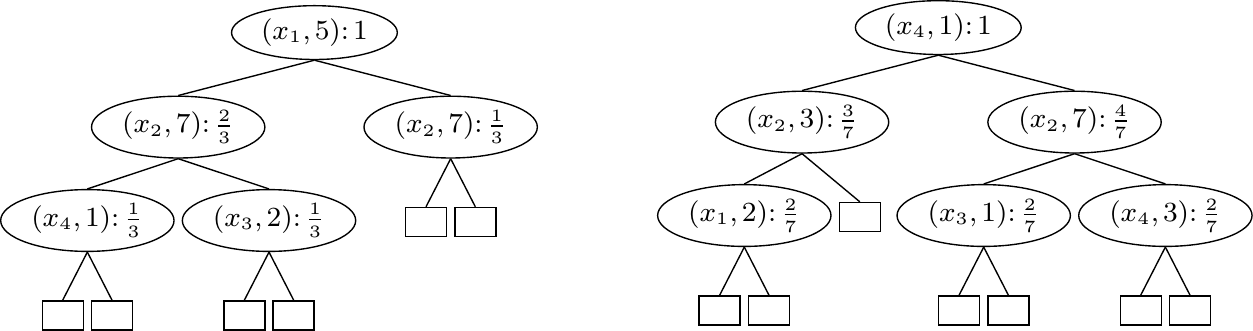}}
	{\Cref{ex:weight_calc} node contributions to \cref{eq:weight_func_cover} weight calculation. Each split node contains `$(x_i, v):w'$ where $(x_i, v)$ is the split pair and $w$ is the node's contribution to $(x_i, v)$'s weight. We calculate $w$ as the proportion of leaves covered relative to the total number of leaves.\label{fig:tree_pseudo_calc}}
	{}
\end{figure}
	\begin{example}\label{ex:weight_calc}
		\Cref{fig:tree_pseudo_calc} shows the weight given to each node for two trees. 
		The left tree contains 6 leaves and the right tree contains 7 leaves.
		Consider split $(x_2, 7)$. The left tree contains two nodes splitting on $(x_2, 7)$ one of which covers 4 out of 6 leaves and the other covers 2 out of 6 leaves therefore these nodes contribute $\frac{2}{3}$ and $\frac{1}{3}$, respectively, to the weight.
		Similarly, the right tree contains a single node splitting on $(x_2, 7)$ which covers 4 out of 7 leaves therefore this node contributes $\frac{4}{7}$ to the weight.
		We obtain the weight for $(x_2,7)$ by summing these values, i.e., $\weight((x_2, 7), \mathcal{T})=\frac{2}{3}+\frac{1}{3}+\frac{4}{7} = 1\frac{4}{7}$.
	\end{example}
	The \cref{eq:weight_func_cover} weight function initializes pseudocosts satisfying the following properties:
	\begin{enumerate}
		\item for each tree $t$, $\weight((x_i, v), t)$ is proportional to $\sum_{s \in r((x_i,v), t)}|\cover(s,t)|$,
		\item if $(x_i, v)$ and $(x_{i'}, v')$ cover the same set of leaves in tree $t$ then $\weight((x_i,v), t) = \weight((x_{i'}, v'), t)$.
	\end{enumerate}
	
\subsubsection{Strong Branching}
\label{subsec:strong_branch}
Branch selection is fundamental to any B\&B algorithm. 
\emph{Strong branching} selects a branch that enables pruning with low effort computations and achieves a non-negligible speed-up in the algorithm's performance
\citep{MORRISON201679}.
Strong branching increases the size of efficiently solvable large-scale mixed-integer problems and is a major solver component \citep{Klabjan2001,Anstreicher2002,Anstreicher2003,10.1007/978-3-540-45157-0_6,doi:10.1080/10556780903087124,misener-floudas:2012,kilinic2014}.
Here, strong branching leverages the easy-to-solve convex penalty term for pruning.

\begin{algorithm}[t]
\caption{Strong Branching}
\label{alg:strong_branch}
\begin{algorithmic}[1]
	\State $S$: B\&B node with bounds $b^{\text{GBT},S}$ and $b^{\text{cvx},S}$
	\State $B^S=[(x_{i_1},v_1),\ldots,(x_{i_l},v_l)]$: $l$ next branches list w.r.t. \cref{subsec:preprocess} pseudo-cost order 
	\For{$(x_i,v)\in B^S$}
	\State $S_{\text{left}},S_{\text{right}}$: $S$ children by branching on $(x_i, v)$
	\State Compute $b^{\text{cvx},S_{\text{left}}}$ and $b^{\text{cvx},S_{\text{right}}}$
	\If {$\max\{b^{\text{cvx},S_{\text{left}}},b^{\text{cvx},S_{\text{right}}}\}+b^{\text{GBT},S}<f^*$}
	\State \Return $\argmin\{b^{\text{cvx},S_{\text{left}}}, b^{\text{cvx},S_{\text{right}}}\},(x_i,v)$
	\EndIf
	\EndFor
	\State \Return $S, (x_{i_1},v_1)$
\end{algorithmic}
\end{algorithm}

\begin{figure}[t]
	\FIGURE
	{\includegraphics[width=0.6\linewidth]{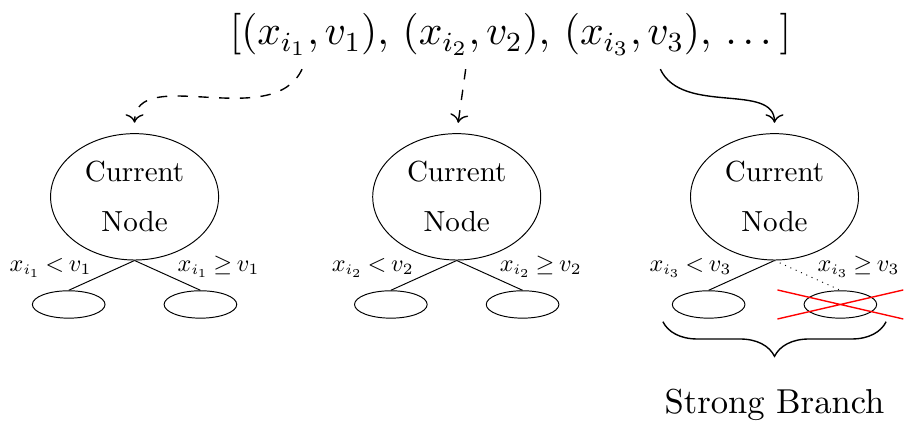}}
	{Strong branching for selecting the next spatial branch. A strong branch leads to a node that is immediately pruned, based on a convex bound computation.\label{fig:strong_branch}}
	{}
\end{figure}

At a B\&B node $S$, branching produces two children $S_{\text{left}}$ and $S_{\text{right}}$.
Strong branching \cref{alg:strong_branch} considers the branches 
in their \cref{subsec:preprocess} pseudo-cost ordering and
assesses each branch by computing the associated convex bound.
Under the strong branching test, one node among $S_{\text{left}}$ and $S_{\text{right}}$ inherits 
the convex bound $b^{\text{cvx},S}$ from the parent, 
while the other requires a new computation.
Suppose that $S'\in\{S_{\text{left}},S_{\text{right}}\}$ does not inherit $b^{\text{cvx},S}$.
If $b^{\text{cvx}, S'}$ satisfies 
the \cref{eq:node_prune} pruning condition without GBT bound improvement, 
then $S'$ is immediately selected as the strong branch and strong branching repeats at the other child node $S''$.
\Cref{fig:strong_branch} illustrates strong branching.
When \cref{alg:bb_overview} does not find a strong branch, 
	it performs a GBT lower bound update and branches on the first item of the branch ordering.
\Cref{alg:bb_overview} then adds this node's children to a set of 
	unexplored nodes and continues with the next B\&B iteration.

Strong branching allows efficient pruning when the convex objective part is significant.
Strong branching may reduce
the computational overhead incurred by GBT bound recalculation when 
\cref{alg:strong_branch} selects multiple strong branches  
between GBT bound updates.
While a single strong branch assessment is negligible, the cumulative cost 
of calculating convex bounds for all branches may be high.
\Cref{subsec:preprocess} orders the branches according to a measure of effectiveness
aiding GBT bounding, so the time spent deriving strong branches with
small weighting function may be better utilized in improving the GBT bound.
Opposed to full strong branching, i.e., assessing all branches, strong branching \cref{alg:strong_branch} uses a lookahead approach \citep{ACHTERBERG200542}.
Parameterized by a lookahead value $l\in\mathbb{Z}_{>0}$, 
\cref{alg:strong_branch} investigates the first $l$ branches. 
If \cref{alg:strong_branch} finds a strong branch,
\cref{alg:bb_overview} repeats \cref{alg:strong_branch}, 
otherwise the B\&B \cref{alg:bb_overview}
updates the GBT bound $b^{\text{GBT},S,P}$ at the current node.
\Cref{alg:strong_branch} keeps strong branching checks relatively cheap 
and maintains a balance between searching and bounding.

\subsection{Heuristics}\label{subsec:feasibility}
To prune, i.e., satisfy \cref{eq:node_prune}, 
consider two heuristic methods generating good feasible solutions 
to Problem~(\ref{eq:full_formulation}):
(i) a mixed-integer convex programming (convex MINLP) approach, and (ii) particle swarm optimization (PSO) \citep{494215,488968}.
The mixed-integer approach uses the decomposability of GBT ensembles, i.e., while convex MINLP solvers provide weak feasible solutions for large-scale instances of Problem~(\ref{eq:full_formulation}), they may efficiently solve moderate instances to global optimality \citep{WESTERLUND1995131,ts:05,vigerske:2012,Misener2014,LUNDELL20172137}.
The PSO approach exploits trade-offs between the convex and objective GBT parts.
Metaheuristics like particle swarm optimization and simulated annealing \citep{Kirkpatrick671} may produce heuristic solutions in preprocessing, i.e., before the branch-and-bound algorithm begins. Simpler convex MINLP heuristics may improve upper bounds at a branch-and-bound node because of their efficient running times.
\Cref{subsec:heuristics} in the electronic companion discusses these heuristics.

\section{Case Studies: Principal Component Analysis for Penalizing Solutions far from Training Data}\label{sec:case_studies}
Our case studies consider GBT instances where training data is not evenly distributed over the $[\gbtVec{v}^L,\gbtVec{v}^U]$ domain.
So, while $\gbtVec{x}\in[\gbtVec{v}^L,\gbtVec{v}^U]$ is feasible, 
$\gbtFunc(\gbtVec{x})$ may be less meaningful for
$\gbtVec{x}$ far from training data. 
The Problem~(\ref{eq:full_formulation}) $\gbtCvx(\gbtVec{x})$ function, for the case studies, is a penalty function constructed with principal component analysis (PCA) \citep{pca}.

PCA characterizes a large, high-dimensional input data set
$D=\{\gbtVec{d}^{(1)}, \dots, \gbtVec{d}^{(p)}\}$
with 
	a low-dimensional subspace capturing most of the variability \citep{introStatLearn}.
PCA defines a set of $n$ ordered, orthogonal \emph{loading vectors}, $\phi_i$, such that $\phi_i$ captures more variability than $\phi_{i'}$, for $i<i'$.
PCA on $D$ defines parameters $\gbtVec{\mu},\gbtVec{\sigma}\in\mathbb{R}^n$ and  $\gbtMat{\Phi}=[\phi_1\,\dots\,\phi_n]\in\mathbb{R}^{n\times n}$, i.e., the sample mean, sample standard deviation and loading vectors, respectively.
Vectors $\gbtVec{\mu}$ and $\gbtVec{\sigma}$ standardize $D$ since PCA is sensitive to scaling.
Often, only a few ($k<n$) leading loading vectors capture most of the variance in $D$ and
$\gbtMat{\Phi}'=[\phi_1\,\dots\,\phi_k]$ may effectively replace $\gbtMat{\Phi}$. $\gbtMat{P}=\gbtMat{\Phi}'\gbtMat{\Phi}^{\prime \top}$ defines a projection matrix to the subspace spanned by $\{\phi_1,\dots,\phi_k\}$.
Penalizing solutions further from training data with PCA defined projection matrix $\gbtMat{P}$:
\begin{equation}\gbtCvx\nolimits_\lambda(\gbtVec{x})=\lambda\left\lVert(\gbtMat{I}-\gbtMat{P})\diag(\gbtVec{\sigma})^{-1}(\gbtVec{x}-\gbtVec{\mu})\right\rVert^2_2\label{eq:cvx_pca}\end{equation}
where $\lambda>0$ is a penalty parameter, $\gbtVec{I}$ is the identity matrix and $\diag(\cdot)$ is a matrix with the argument on the diagonal. 
Larger $\lambda$ is more conservative with respect to PCA subspace $\gbtMat{P}$.
Note in \cref{eq:cvx_pca} that our specific nonlinear convex penalty is a convex quadratic.

	\Cref{eq:cvx_pca} aims to characterize the region containing the training data with an affine subspace.
	Points in the subspace are not penalized and points close to the subspace are not heavily penalized.
	However, \cref{eq:cvx_pca} may be qualitatively less effective when the standardized training data is not evenly distributed within subspace $\gbtMat{P}$.
	\begin{example}\label{eg:cluster}
		Consider a data set $\left\{\bm{x}^{(i)}\right\}_{i=1}^{2m}$, $\bm{x}^{(i)}\in\mathbb{R}^3$ where  $x_1^{(i)}\sim U(0,1), \, \forall \ i\in[2m]$, and $x_2^{(i)}=x_3^{(i)}=0$, $x_2^{(m+i)}=x_3^{(m+i)}=1, \, \forall \ i\in[m]$. The 2D subspace containing these points contains the origin and directions $(1,0,0)^T, (0, 1, 1)^T$. 
	\Cref{eq:cvx_pca} does not penalize points in this subspace.
		But the point $(0.5, 0.5, 0.5)$, which is contained in the subspace, is far from the training data when considering the subspace distribution.
		Having $x_2^{(i)},x_3^{(i)}\sim N(0,\varepsilon)$, $x_2^{(m+i)},x_3^{(m+i)}\sim N(1,\varepsilon), \, \forall \ i\in[m]$ and small $\varepsilon>0$, introduces an error term to the second and third variables while retaining the same clustered distribution over the subspace.
	\end{example}
	Clustering, e.g., \cref{eg:cluster}, may be handled by the \cref{sec:branch_and_bound} B\&B. We could instantiate a separate instance for each cluster using a penalty that only considers training data in a given cluster and limit the solve to a reduced box domain. 
		A single problem formulation considering more complex training data relationships may 
negatively affect the strong branching aspect of B\&B \cref{alg:bb_overview}.

\section{Numerical Results}
\label{sec:numerical_analysis}

This section compares the \cref{sec:branch_and_bound} lower bounding and branch-and-bound algorithms to black-box solvers.
\Cref{app:feasibility} of the electronic companion presents results for the \cref{subsec:feasibility} heuristics.
\Cref{Section:Specifications} provides information about the system specifications and the solvers.
\Cref{subsec:concrete,subsec:basf} investigate two GBT instances for engineering applications, namely: (i) concrete mixture design and (ii) chemical catalysis.
\Cref{subsec:observations} discusses observations from the \cref{subsec:concrete,subsec:basf} results.
The concrete mixture design instance is from the UCI machine learning repository \citep{Dua:2017}. 
The industrial chemical catalysis instance is provided from BASF.
\Cref{tab:bb_parameters} presents information about these instances.
For both instances, we model closeness to training data using the PCA-based function $\gbtCvx(\bm{x})$ defined in \cref{eq:cvx_pca}.

\begin{table}
	\TABLE
	{Instance Sizes \label{tab:bb_parameters}}	
	{\begin{tabular}{llcccc}
	\toprule
	& && Concrete Mixture Design & Chemical Catalysis\\
	\midrule
	\multicolumn{2}{l}{\textit{GBT attributes:}}\\
	&Number of trees && 7,750 & 8,800\\
	&Maximum depth && 16 & 16\\
	&Number of leaves && 131,750 & 93,200\\
	&Number of $x_i$ continuous variables && 8&42\\[0.8ex]
	\multicolumn{2}{l}{\textit{Convex MINLP (\ref{Eq:Optimization_Problem}) attributes:}}\\
	&Number of $y_{i,j}$ binary variables && 8,441&2,061\\
	&Number of constraints && 281,073&183,791\\[0.8ex]
	\bottomrule
\end{tabular}
}
	{}
\end{table}

\subsection{System and Solver Specifications}
\label{Section:Specifications}

Experiments are run on an Ubuntu 16.04 HP EliteDesk 800 G1 TWR with 16GB RAM and an Intel Core i7-4770@3.40GHz CPU.
Implementations are in Python 3.5.3 using Pyomo 5.2 \citep{hart2011pyomo,hart2017pyomo} for mixed-integer programming modeling and interfacing with solvers.
We use CPLEX 12.7 and Gurobi 7.5.2 as: (i) black-box solvers for the entire convex MINLP (\ref{Eq:Optimization_Problem}), and
(ii) branch-and-bound algorithm components for solving MILP (\ref{eq:ilp_sub}) instances in the \cref{subsec:bb_lb} GBT lower bounding procedure.
Note that current versions of CPLEX and Gurobi cannot solve general convex MINLP, so we would use a more general solver if we had non-quadratic penalty functions.
All results report wall clock times.

This section evaluates the (i) objective lower bounding procedure, and (ii) branch-and-bound algorithm, both of which use CPLEX or Gurobi as a black-box MILP solver. 
We also apply CPLEX and Gurobi to the entire MINLP for evaluating branch-and-bound \cref{alg:bb_overview}. 
Figures 5-12 append labels -C and -G to indicate CPLEX and Gurobi, respectively, and use different line types for displaying the results.
At nodes immediately following a GBT bound update, the B\&B algorithm assesses solutions from solving the convex part of Problem~(\ref{eq:full_formulation}) as heuristics solutions. 
We use the default CPLEX 12.7 and Gurobi 7.5.2 tolerances, i.e., relative MIP gap, integrality and barrier convergence tolerances of $10^{-4}$, $10^{-5}$ and $10^{-8}$, respectively. 

\subsection{Concrete Mixture Design}\label{subsec:concrete}
In concrete mixture design, different ingredient proportions result in different 
	properties of the concrete, e.g., compressive strength.
The relationship between ingredients and properties is complex, so black-box machine learning is 
	well suited for the function estimation task \citep{doi:10.1061/(asce)cp.1943-5487.0000088,erdal20131689,DEROUSSEAU201842}.

\subsubsection{Instance}
We maximize concrete compressive strength where GBTs are used for 
modeling.
Since we maximize concrete compressive strength, negating all leaf weights $F_{t,l}$ forms  
	an equivalent GBT instance that fits the Problem~(\ref{eq:full_formulation}) minimization formulation.
We use the \citet{YEH19981797} concrete compressive strength dataset from the 
	UCI machine learning repository \citep{Dua:2017}.
This dataset has $n=8$ continuous variables.
R packages gbm \citep{Ridgeway} and caret \citep{JSSv028i05} are used for GBT training.
Root-mean-square error is used for model selection.
The resulting GBT instance has 7750 trees with max depth 16.
The PCA based convex penalty has $\gbtRank(\gbtMat{P})=4$, i.e., we select the first four loading vectors.
\Cref{subsubsec:concrete_heur} of the electronic companion presents the lesults for the \cref{subsec:feasibility} heuristics.

\subsubsection{GBT Lower Bounding}\label{subsubsec:concrete_bounding}
\begin{figure}
	\centering
	\includegraphics[width=\linewidth,height=4cm]{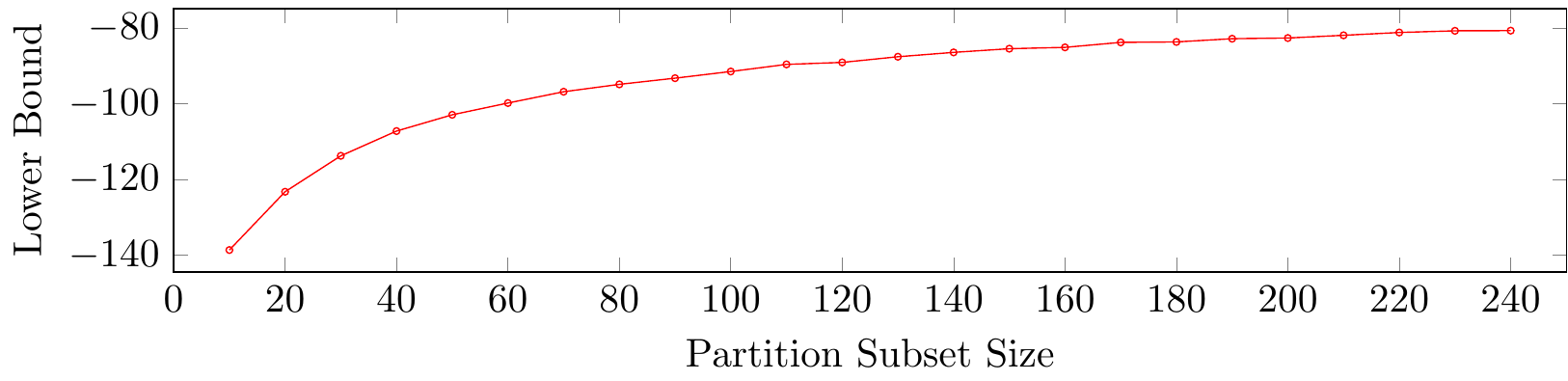}
	\caption{%
		Concrete mixture design instance: Global GBT lower bound improvement using the \cref{subsec:bb_lb} GBT lower bounding approach for different partition subset sizes.
	\label{fig:concrete_grp_bnd_bnd}}
\end{figure}
\begin{figure}
	\centering
	\includegraphics[width=\linewidth,height=6cm]{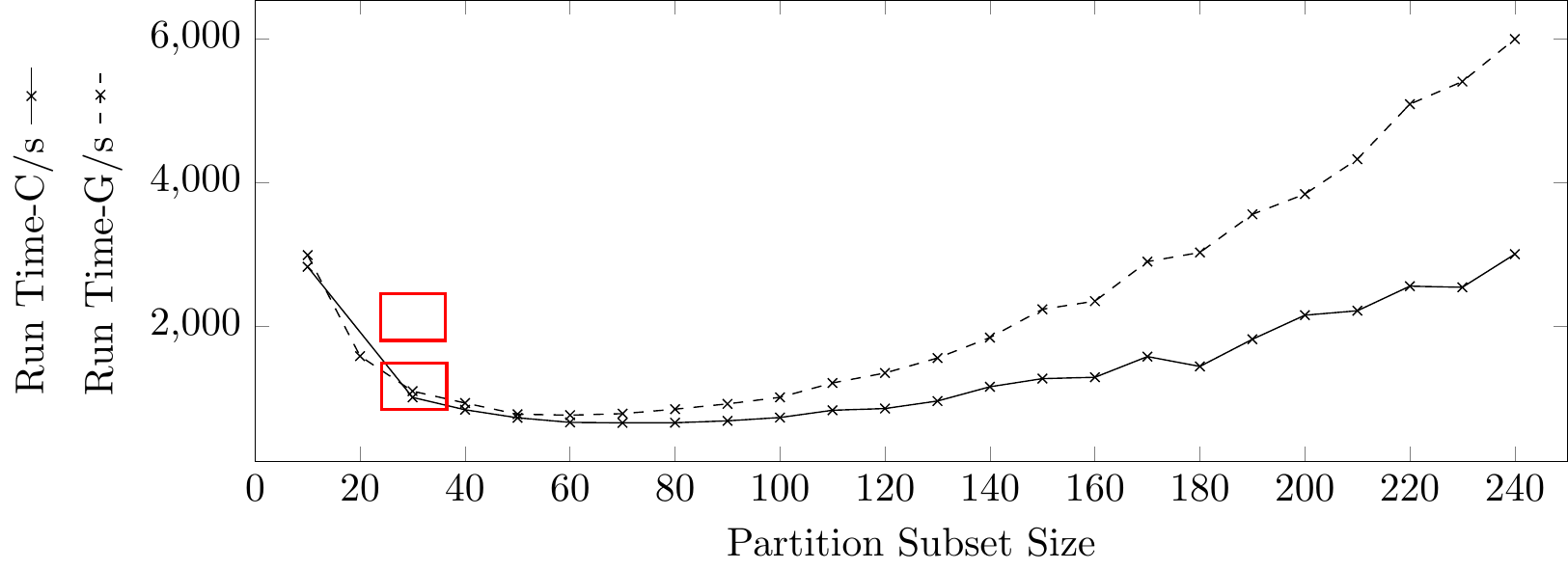}
	\caption{%
		Concrete mixture design instance: Global GBT lower bounding wall clock time using the \cref{subsec:bb_lb} approach for different partition subset sizes.
		Suffixes -C and -G denote subsolvers CPLEX 12.7 and Gurobi 7.5.2, respectively.
	\label{fig:concrete_grp_bnd_time}}
\end{figure}
\Cref{fig:concrete_grp_bnd_bnd,fig:concrete_grp_bnd_time} evaluate the \cref{subsubsec:gbt_bound} GBT lower bounding approach for different partition subset sizes. 
\Cref{fig:concrete_grp_bnd_bnd} illustrates the global GBT lower bound improvement as the partition subset size increases.
\Cref{fig:concrete_grp_bnd_time} compares run times with either CPLEX 12.7, or Gurobi 7.5.2 as subsolvers for each partition subset size.
For the entire MILP instance, i.e.,\ solving Problem~(\ref{eq:ilp_sub}), black-box solving with CPLEX 12.7 and Gurobi 7.5.2 achieve GBT lower bounds -97 and -547, respectively, within 1 hour.
The \cref{subsubsec:gbt_bound} approach achieves a lower bound of -83 (partition size 190), in 1 hour, and improves upon black-box solver lower bounds in under 15 minutes (partition size 70).

\subsubsection{Branch-and-Bound Algorithm}\label{subsubsec:concrete_bb}
\begin{figure}
	\centering
	\includegraphics[width=0.98\linewidth, height=24ex]{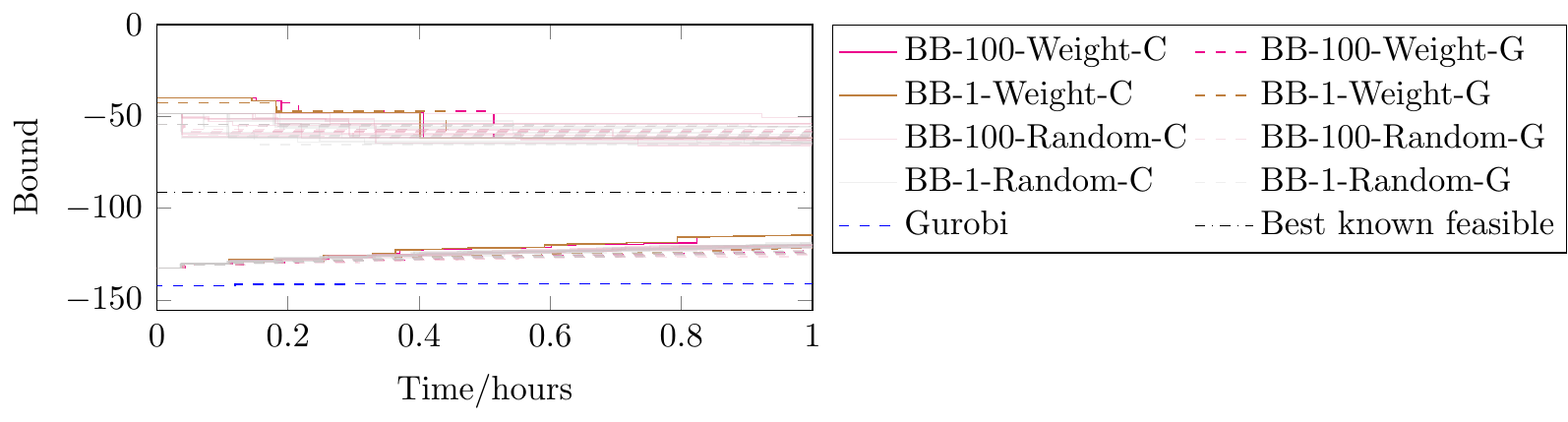}
	\caption{%
			Concrete mixture design instance ($\lambda=1$): B\&B lower bound improvement compared to Gurobi 7.5.2, with a one hour timeout.
	The B\&B \cref{alg:bb_overview} is labeled BB-$a$-$b$-$c$ where $a$, $b$ and $c$  denote the strong branching lookahead value, the pseudocost initialization approach, and the solver used for lower bounding and solving convex quadratics, respectively. The BB-\textasteriskcentered{} results sort the unexplored nodes in ascending lower bound order. The dashed-dotted line reports best found feasible solution (upper bound).
	\label{fig:concrete_bnd_evolution_lm_1}}
\end{figure}

\begin{figure}
	\centering
	\includegraphics[width=0.98\linewidth, height=24ex]{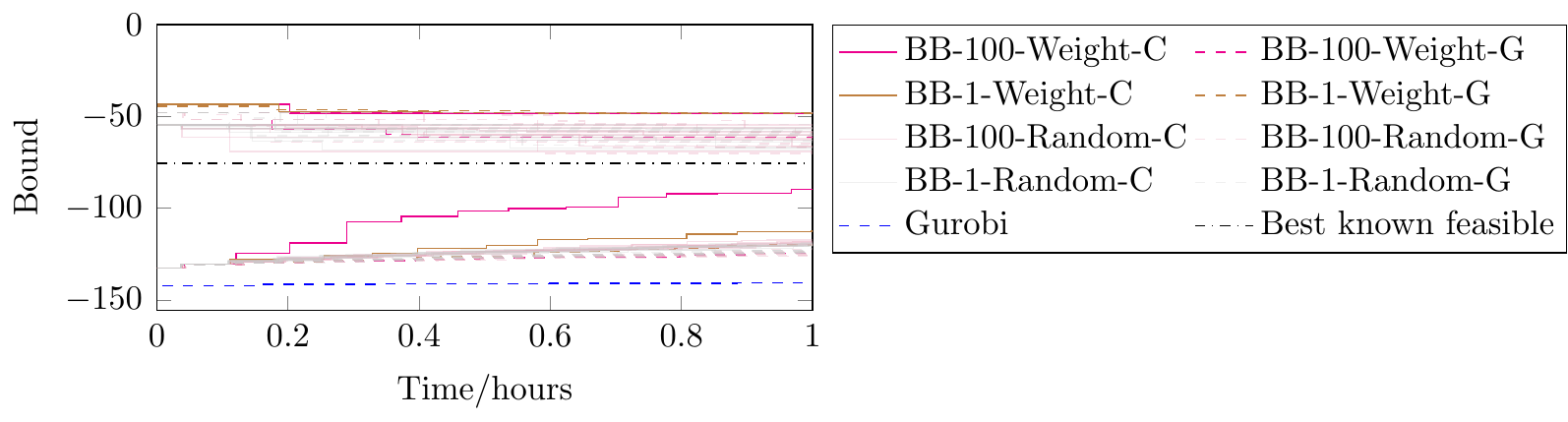}
	\caption{%
			Concrete mixture design instance ($\lambda=1000$): B\&B lower bound improvement compared to Gurobi 7.5.2, with a one hour timeout.
	The B\&B \cref{alg:bb_overview} is labeled BB-$a$-$b$-$c$ where $a$, $b$ and $c$  denote the strong branching lookahead value, the pseudocost initialization approach, and the solver used for lower bounding and solving convex quadratics, respectively. The BB-\textasteriskcentered{} results sort the unexplored nodes in ascending lower bound order. The dashed-dotted line reports best found feasible solution (upper bound).
	\label{fig:concrete_bnd_evolution_lm_1000}}
\end{figure}

	We instantiate the branch-and-bound algorithm with a root node partition of 70 trees, and non-root lower bounding time limit of 120 seconds. All branch-and-bound tests are run with CPLEX 12.7 and Gurobi 7.5.2 as MILP subsolvers. We assess the effect of strong branching by comparing lookahead list sizes $l = 1$ vs.\ $l = 100$. We assess the quality of feasible solutions by comparing with the \cref{tab:concrete_feas_results} [electronic companion] best found feasible solution. We assess the pseudocost ordering by comparing with 10 independent tests of random branch orderings for each strong branch lookahead-subsolver combination. We compare all branch-and-bound results, which allocate 1 hour for GBT lower bounding at the root node and 1 hour for the B\&B search, to 3 hour black-box runs of CPLEX 12.7 and Gurobi 7.5.2 for the entire convex MINLP.

	\Cref{fig:concrete_bnd_evolution_lm_1,fig:concrete_bnd_evolution_lm_1000} plot the bound improvement for $\lambda=1$ and $\lambda=1000$, respectively. For the entire convex MINLP, the black-box CPLEX 12.7 bounds are outside the figure axis limits. For $\lambda=1$, a larger strong branching lookahead value does not noticeably improve the lower bound, but a larger lookahead does significantly improve the lower bound for $\lambda=1000$.
\Cref{fig:concrete_bnd_evolution_lm_1} depicts the lower bound improvement. 
The B\&B algorithm lower bound improves over time, but there is still a non-negligible gap from the best-known feasible solution after 1 hour.
This gap appears to be due to a cluster-like effect caused by the GBTs \citep{Du1994,Wechsung2014,Kannan2017}, where the variable split points are quite close. 
In the B\&B algorithm, if the current lookahead list contains these clusters, strong branching is less effective.
CPLEX 12.7 results in an out-of-memory error prior to beginning the branch-and-bound search therefore its lower bounds are relatively poor.
Gurobi 7.5.2 returns an incumbent of -85 and a lower bound of -141, after 2 hours, and these do not improve further in the subsequent hour. 
The B\&B algorithm, at 2 hours, i.e., prior to tree search, has an incumbent of -91 and a lower bound not less than -133. 
Given an additional hour for tree search, the gap reduces further.
\Cref{tab:concrete_24hr} compares the B\&B algorithm to Gurobi 7.5.2 with 24 hours time limit. 
The Gurobi heuristics generally outperform the B\&B algorithm, but the B\&B algorithm derives better lower bounds.
In all cases, $\geq 22$\% optimality gap remains. 
Because regions close to training data have many GBT breakpoints, optimal solutions lie in highly discretized areas of the feasibility domain.

\begin{table}
	\centering
	\TABLE
	{Concrete mixture design instance: Results comparing 24 hour runs of the B\&B algorithm with Gurobi 7.5.2. The B\&B algorithm uses a strong branching lookahead value of 100, a root node partition of 70 trees, a non-root lower bounding time limit of 120 seconds, and CPLEX 12.7 as a subsolver.\label{tab:concrete_24hr}}
	{\begin{tabular}{lclrrrlrrrl}
	\toprule
	&&&\multicolumn{3}{c}{BB-C}&&\multicolumn{3}{c}{Gurobi 7.5.2}&\\
	&$\lambda$ && UB & LB &Gap &&UB&LB & Gap\\
	\midrule
	&$1$    &&  $-80.56$ & $-99.87$ & $24\%$ && $-85.48$ & $-140.75$&$64\%$\\
	&$10$   &&  $-74.96$ & $-99.39$ & $33\%$ && $-85.06$ & $-121.10$&$42\%$\\
	&$100$  &&  $-73.74$ & $-96.43$ & $31\%$ && $-77.98$ & $-121.27$&$55\%$\\
	&$1000$ &&  $-74.86$ & $-90.75$ & $22\%$&& $-72.29$ & $-121.23$&$67\%$\\
	\bottomrule
\end{tabular}
}
	{}
\end{table}

\subsection{Chemical Catalysis}\label{subsec:basf}
BASF uses catalysts to improve yield and operating efficiency.
But, modeling catalyst effectiveness is highly nonlinear 
	and varies across different applications.
BASF has found GBTs effective for modeling catalyst behavior.
Capturing the high-dimensional nature of catalysis over the entire feasible domain requires many experiments, too many to run in practice.
Running a fewer number of experiments necessitates penalizing solutions further from where the GBT function is trained.

\subsubsection{Instance}
The BASF industrial instance contains $n=42$ continuous variables.
The convex part of the instance takes the following form:
\begin{equation}
	\gbtCvx\nolimits_\lambda(\gbtVec{x})=\lambda\left\lVert(\gbtMat{I}-\gbtMat{P})\diag(\gbtVec{\sigma})^{-1}(\gbtVec{x}-\gbtVec{\mu})\right\rVert^2_2 + \left(100-\sum_{i\in\mathcal{I}^{\%}}x_i\right)^2\label{eq:basf_cvx}
\end{equation}
\Cref{eq:basf_cvx} differs from \cref{eq:cvx_pca} in its addend which aims to generate solutions where $x_i\in\mathcal{I}^\%$, i.e., proportions of 
	the chemicals being mixed, sum to $100\%$.
The test instance has $\gbtRank(\gbtMat{P})=2$ and $|\mathcal{I}^{\%}|=37$.
The GBT part contains 8800 trees where 4100 trees have max depth 16,
	the remaining trees have max depth 4, the total number of leaves is 93,200 and the corresponding 
	Problem~(\ref{eq:ilp_sub}) MILP model has 2061 binary variables.
\Cref{subsubsec:basf_heur} of the electronic companion presents the results for the \cref{subsec:feasibility} heuristics.

\subsubsection{GBT Lower Bounding}\label{subsubsec:basf_bounding}
\begin{figure}
	\FIGURE
	{\includegraphics[width=\linewidth,height=4cm]{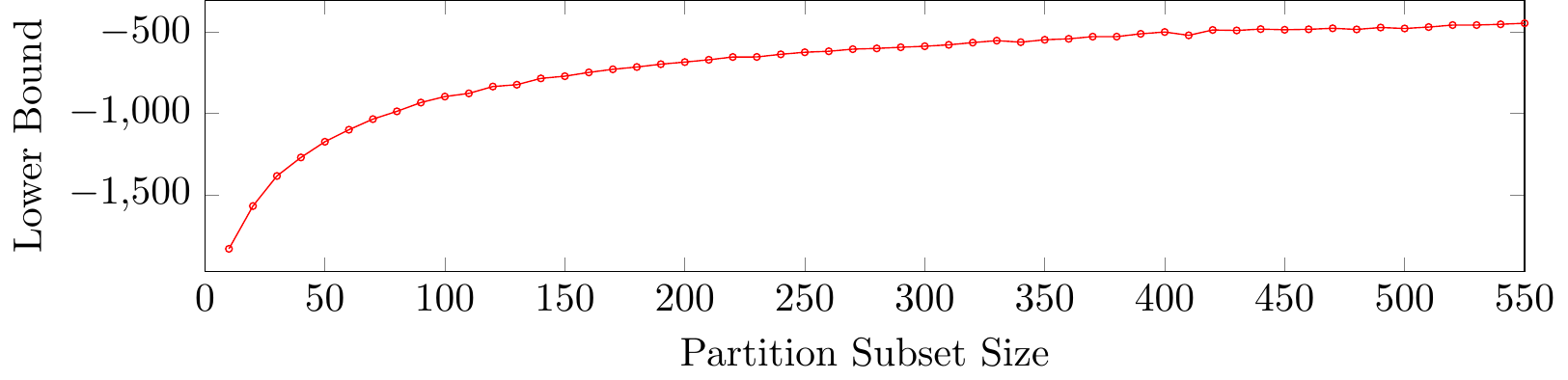}}
	{
		Chemical catalysis BASF instance: Global GBT lower bound improvement using the \cref{subsec:bb_lb} GBT lower bounding approach  for different partition subset sizes. 
	\label{fig:basf_standard_grp_bnd_bnd}
	}
 {}
\end{figure}

\begin{figure}
	\FIGURE
	{\includegraphics[width=\linewidth,height=5cm]{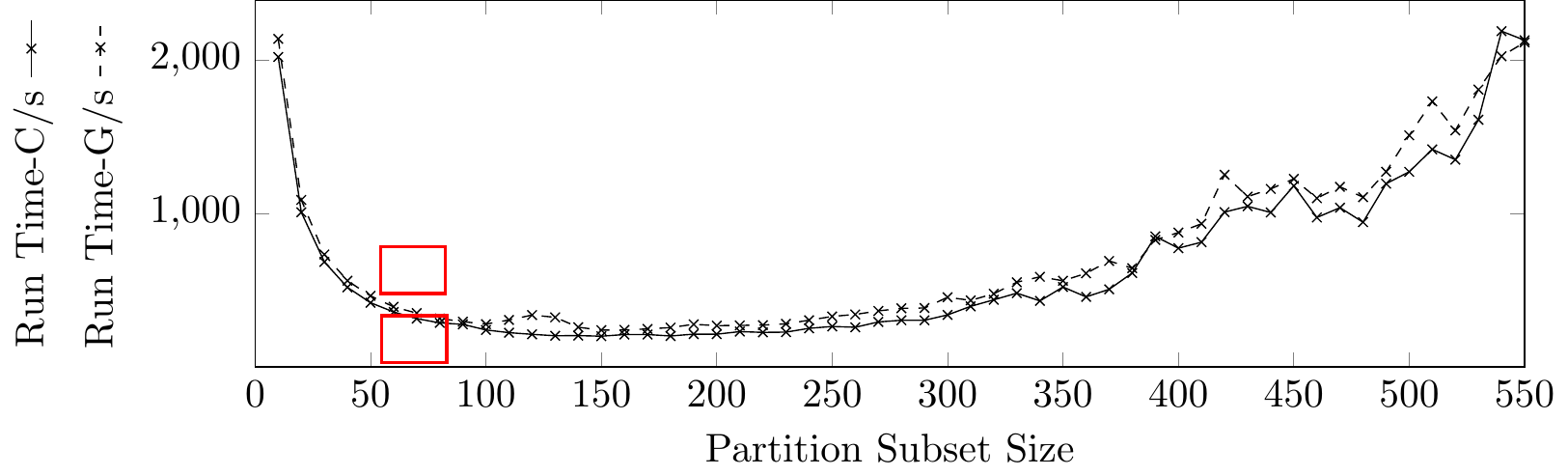}}
	{
		Chemical catalysis BASF instance: Global GBT lower bounding wall clock time using the \cref{subsec:bb_lb} approach for different partition subset sizes. 
		Suffixes -C and -G denote subsolvers CPLEX 12.7 and Gurobi 7.5.2, respectively.
	\label{fig:basf_standard_grp_bnd_time}
	}
 {}
\end{figure}

\Cref{fig:basf_standard_grp_bnd_bnd,fig:basf_standard_grp_bnd_time} evaluate the \cref{subsubsec:gbt_bound} GBT lower bounding approach for different partition subset sizes.
\Cref{fig:basf_standard_grp_bnd_bnd} illustrates the global GBT lower bound improvement as the partition subset size increases.
\Cref{fig:basf_standard_grp_bnd_time} compares run times when using either CPLEX 12.7, or Gurobi 7.5.2 as subsolvers for each partition subset size.
These results resemble \cref{fig:concrete_grp_bnd_bnd,fig:concrete_grp_bnd_time}.
In particular, (i) the lower bound is improved with larger subset sizes, (ii) there is a time-consuming modeling overhead for solving many small MILPs for small subset sizes, and (iii) the running time increases exponentially, though non-monotonically, for larger subset sizes.
We compare the lower bounding approach with solving the entire MILP (\ref{eq:ilp_sub}) using CPLEX 12.7, or Gurobi 7.5.2 as black-box solvers.
Our lower bounding approach exhibits a superior time-to-lower bound performance:
(i) it improves the Gurobi 7.5.2 lower bound with subset size 140 and 4 minutes of execution, and (ii) it improves the CPLEX 12.7 lower bound with subset size 360 and 8 minutes of execution.

\subsubsection{Branch-and-Bound Algorithm}\label{subsubsec:basf_bb}
We instantiate the branch-and-bound algorithm with a root node partition of 150 trees, and non-root lower bounding time limit of 120 seconds. All branch-and-bound tests are run with CPLEX 12.7 and Gurobi 7.5.2 as subsolvers. We assess the effect of strong branching by comparing lookahead list sizes $l = 1$ vs.\ $l = 100$. We assess the quality of feasible solutions by comparing with the \cref{tab:basf_feas_results} [electronic companion] best found feasible solution. We assess the pseudocost ordering by comparing with 10 independent tests of random branch orderings for each strong branch lookahead-subsolver combination. 
We compare all branch-and-bound results, which allocate 1 hour for GBT lower bounding at the root node and 1 hour for the B\&B search, to 3 hour black-box runs of CPLEX 12.7 and Gurobi 7.5.2 for the entire convex MINLP.
	
\begin{figure}
	\centering
	\FIGURE
	{
		\includegraphics[height=24ex, width=0.98\linewidth]{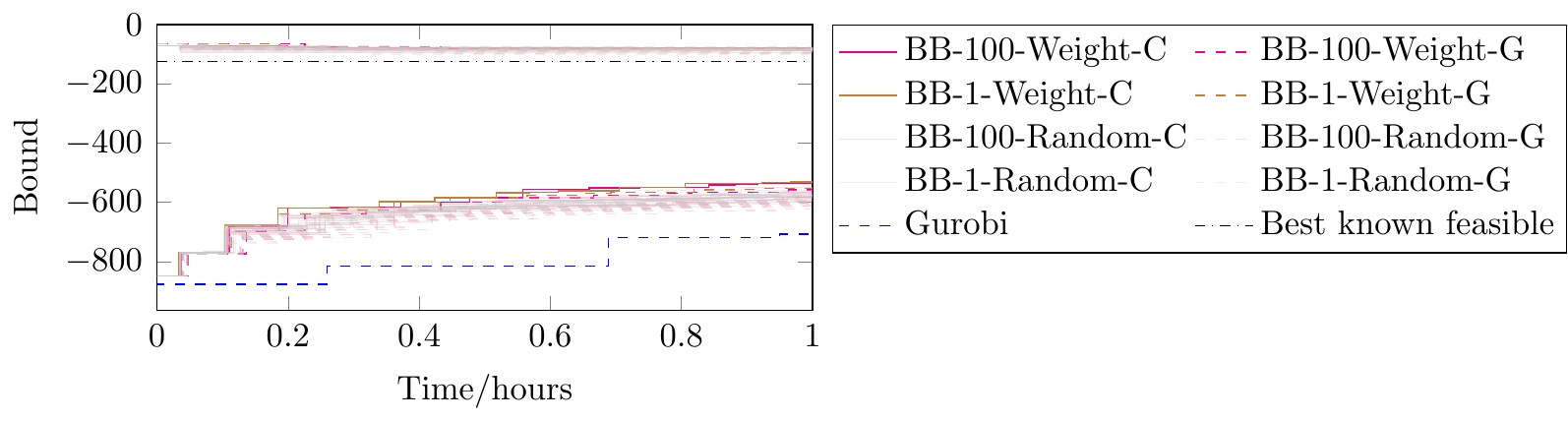}
	}
	{%
Chemical catalysis BASF instance ($\lambda=1$): B\&B lower bound improvement compared to Gurobi 7.5.2 with one hour timeout.
	The B\&B \cref{alg:bb_overview} is labeled BB-$a$-$b$-$c$ where $a$, $b$ and $c$  denote the strong branching lookahead value, the pseudocost initialization approach, and the solver used for lower bounding and solving convex quadratics, respectively. The BB-\textasteriskcentered{} results sort the unexplored nodes in ascending lower bound order. The dashed-dotted line reports best found feasible solution (upper bound).
	\label{fig:bound_evo_lm_1}}
	{}
\end{figure}

\begin{figure}
	\centering
	\FIGURE
	{
		\includegraphics[height=24ex, width=0.98\linewidth]{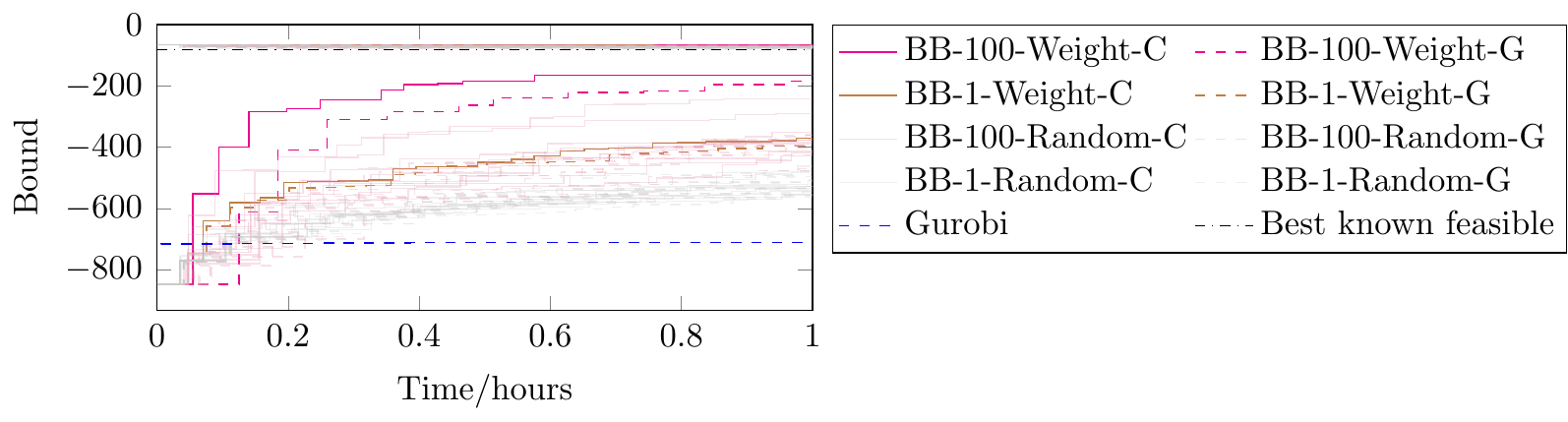}
	}
	{%
		Chemical catalysis BASF instance ($\lambda=1000$): B\&B lower bound improvement compared to Gurobi 7.5.2 with one hour timeout.
	The B\&B \cref{alg:bb_overview} is labeled BB-$a$-$b$-$c$ where $a$, $b$ and $c$  denote the strong branching lookahead value, the pseudocost initialization approach, and the solver used for lower bounding and solving convex quadratics, respectively. The BB-\textasteriskcentered{} results sort the unexplored nodes in ascending lower bound order. The dashed-dotted line reports best found feasible solution (upper bound).
	\label{fig:bound_evo_lm_1000}}
	{}
\end{figure}

\Cref{fig:bound_evo_lm_1,fig:bound_evo_lm_1000} plot the bound improvement
	for $\lambda=1$ and $\lambda=1000$, respectively.
For the entire convex MINLP, CPLEX 12.7 reports a poor lower bound and does not find a feasible solution within 3 hours.
The B\&B algorithm terminates with a tighter lower bound and closes a larger gap than the black-box solvers, across all tested parameter combinations.
The B\&B algorithm performs better for $\lambda=1000$ because the convex part dominates the GBT part more, making strong branching more effective.
Finally, we see that the branch-and-bound algorithm finds a relatively good heuristic solution at the root node for $\lambda=1000$. For $\lambda=1$, there is a larger gap between the B\&B upper bounds and the best known feasible solution this is expected as solving this problem is closer to optimizing only over the GBT MILP where an optimal solution may be further from the PCA subspace.
\Cref{tab:basf_24hr} compares the B\&B algorithm to Gurobi 7.5.2 with 24 hours time limit.
The Gurobi heuristic solutions generally outperform the B\&B algorithm.
Nevertheless, the B\&B algorithm derives better lower bounds.
For $\lambda=1000$, the B\&B algorithm succeeds in proving global optimality whereas Gurobi terminates with a $7\%$ gap.

\begin{table}
	\centering
	\TABLE
	{Chemical catalysis instance: Results comparing 24 hour runs of the B\&B algorithm with Gurobi 7.5.2. The B\&B algorithm uses a strong branching lookahead value of 100, a root node partition of 150 trees, a non-root lower bounding time limit of 120 seconds, and CPLEX 12.7 as a subsolver.\label{tab:basf_24hr}}
	{\begin{tabular}{clrrrlrrr}
	\toprule
	&&\multicolumn{3}{c}{BB-C}&&\multicolumn{3}{c}{Gurobi 7.5.2}\\
	$\lambda$ && Upper Bound & Lower Bound &Gap&&Upper Bound&Lower Bound&Gap\\
	\midrule
	$1$    &&  $-81.7$ & $-366.0$ & $348\%$ && $-154.8$ & $-580.6$ & $275\%$\\
	$10$   &&  $-80.6$ & $-336.8$ & $318\%$ && $-118.8$ & $-577.1$ & $386\%$\\
	$100$  &&  $-87.3$ & $-187.0$ & $114\%$ &&  $-94.2$ & $-424.5$ & $350\%$\\
	$1000$ &&  $-86.0$ &  $-86.0$ & $0\%$ &&  $-85.9$ &  $-92.1$ & $7\%$\\
	\bottomrule
\end{tabular}
}
	{}
\end{table}

\subsection{Observations}\label{subsec:observations}
The \cref{subsubsec:concrete_bounding,subsubsec:basf_bounding} GBT lower bounding results show that, for large-scale GBT instances, selecting an appropriate partition subset size in the decomposition approach results has a better time-to-lower bound performance than 1 hour black-box MILP solvers.
Both problem instances show that, for larger subset sizes, the running time exponentially increases, while the lower bound improvement rate exponentially decreases.
This is an expected result for GBT instances with deep trees as deeper tree induce more infeasible combinations of branches.
For shallower GBT instances, individual trees to may interact less with each other, hence the decomposition strategy may derive a poorer bound than a black-box MILP solver.
For small subset sizes, the partition-based lower bounding has decreasing running time because of the overhead from many sequential subproblems.

The \cref{subsubsec:concrete_bb,subsubsec:basf_bb} B\&B results also show common features.
Comparing the BB-\textasteriskcentered-$b$-\textasteriskcentered{} results for $b\in\{\text{Weight}, \text{Random}\}$ in \cref{fig:concrete_bnd_evolution_lm_1,fig:concrete_bnd_evolution_lm_1000,fig:bound_evo_lm_1,fig:bound_evo_lm_1000} assess the pseudocost effect.
The \cref{eq:weight_func_cover} initialization outperforms random ordering (for matching lookahead values), showing that the pseudocosts select branches that aid GBT lower bounding. This pseudocost effect is more pronounced with a lookahead value of 100 since multiple branches are selected between branch-and-bound iterations.
For $\lambda=1000$, a lookahead list size $l = 100$ closes more gap than $l = 1$ (comparing BB-100-\textasteriskcentered{} to BB-1-\textasteriskcentered{}), as the B\&B algorithm accepts more branches for strong branching. The difference between $l = 100$ and $l = 1$ implies that increased strong branching improves the GBT lower bound earlier and more often. 
For $\lambda=1$, using a larger strong branching lookahead size does not have a noticable effect.
However, this last finding does not depreciate strong branching.
Since the GBT part dominates the convex aspect for small $\lambda$ values, tighter GBT lower bounds might be essential for taking full advantage of strong branching.
Testing the B\&B algorithm and Gurobi 7.5.2 with a 24 hour run time shows that the branch-and-bound algorithm tends to result in superior lower bounds and closes a larger proportion of the optimality gap whereas Gurobi 7.5.2 produces better heuristic solutions.
Closing any outstanding gap proves difficult as the domains of the remaining unexplored nodes are highly discretized by the GBTs.

\section{Discussion}\label{sec:discussion}
	Our optimization problem consists of: (i) the GBTs, and (ii) the PCA-based penalty.
	Functions obtained from limited, known evaluations with machine learning are approximate by default and may deviate from the ground truth, thus, resulting in false optima.
	The final solution error depends on the training data distribution, noise, and machine learning model.
	Our PCA-based approach may deteriorate for clustered data, e.g., \cref{eg:cluster}, when regions of the PCA subspace are far from training observations.
	A remedy is using data analysis, e.g., clustering \citep{elements}, to assess uniformity in the training data distribution.
	An alternative direction is using proximity measures \citep{rnews,2017arXiv170510883M}.
	The proximity measures may require adjusting when using GBTs since the boosting procedure results in some trees being more relevant than others.
Finally, other convex penalties are relevant in a variety of applications \citep{doi:10.1002/aic.690320408}.
		
	Finally, we acknowledge other approaches for decision-making with optimization problems whose input is specified by machine learning models.
	\citet{NIPS2017_7132} consider end-to-end task-based learning where probabilistic models are trained to be subsequently used within stochastic programming tasks.
	\citet{2017arXiv171008005E} develop a framework for training predictive models with a specific loss function so that the resulting optimization problem has desirable convexity properties and is statistically consistent.
	\citet{2018arXiv180905504W} propose a two-stage approach for integrating machine learning predictions with combinatorial optimization problem decisions.
	The main difference with our work is that we are more focused on the optimization side.

\section{Conclusion}\label{sec:conclusion}
As machine learning methods mature,
decision makers want to move from solely 
making predictions on model inputs to 
integrating pre-trained machine learning models into larger decision-making problems.
This paper addresses a large-scale, industrially-relevant gradient-boosted tree model by directly
exploiting: (i) advanced mixed-integer programming technology with strong optimization formulations,
(ii) GBT tree structure with priority towards searching on commonly-occurring variable splits,
and (iii) convex penalty terms with enabling fewer mixed-integer optimization updates.
The general form of the optimization problem appears whenever we wish to optimize a pre-trained gradient-boosted tree with convex terms in the objective, e.g., penalties.
It would have been alternatively possible to train and then optimize a smooth and continuous machine learning
model, but applications with legacy code may start with a GBT.
Our numerical results test against concrete mixture design and chemical catalysis, two applications
	where the global solution to an optimization problem is often particularly useful.
Our methods 
not only generate good feasible solutions to the optimization problem,
	but they also converge towards proving the exact solution.

\ACKNOWLEDGMENT{%
The support of: BASF SE, the EPSRC Centre for Doctoral Training in High Performance Embedded and Distributed Systems to M.M. (EP/L016796/1), and an EPSRC Research Fellowship to R.M. (EP/P016871/1).%
}

\bibliographystyle{informs2014}
\bibliography{citations}

\newpage

\clearpage
\pagenumbering{arabic}
\renewcommand*{\thepage}{A\arabic{page}}

\begin{APPENDICES}
\crefalias{section}{appsec}
\begin{center}

\begin{large}
Electronic supplementary material: Mixed-Integer Convex Nonlinear Optimization with Gradient-Boosted Trees Embedded \\[10pt]
\end{large}

Miten Mistry \textbullet{} Dimitrios Letsios \textbullet{} Gerhard Krennrich \textbullet{} Robert M.\ Lee \textbullet{} Ruth Misener \\[10pt]
\end{center}

\section{Full convex MINLP formulation}\label{app:full_micp}

\begin{subequations}
\begin{align}
	\min_{\gbtVec{v}^L\leq\gbtVec{x}\leq\gbtVec{v}^U}\;&\gbtCvx(\gbtVec{x}) + \sum_{t\in\mathcal{T}}\sum_{l\in\mathcal{L}_t}F_{t,l}z_{t,l}\\
	\text{s.t.}\;
		&\sum_{l\in\mathcal{L}_{t}}z_{t,l} = 1, &\forall t&\in \mathcal{T},\\
		&\sum_{\mathclap{l\in\gbtLeft_{t,s}}}z_{t,l} \leq y_{i(s),j(s)}, &\forall t&\in\mathcal{T}, s\in 
		\mathcal{V}_{t},\\
		&\sum_{\mathclap{l\in\gbtRight_{t,s}}}z_{t,l} \leq 1-y_{i(s),j(s)}, &\forall t&\in\mathcal{T}, s\in \mathcal{V}_{t},\\
		&y_{i,j}\leq y_{i,j+1}, &\forall i&\in[n],\,j\in[m_i-1],\\
		&x_i\geq v_{i,0} + \sum_{j=1}^{m_i}(v_{i,j} - v_{i,j-1})(1-y_{i,j}),&\forall i&\in[n],\\
		&x_i\leq v_{i,m_i+1}+ \sum_{j=1}^{m_i}(v_{i,j} - v_{i,j+1})y_{i,j},&\forall i&\in[n],\\
		&y_{i,j}\in\{0,1\}, &\forall i&\in[n],\, j\in [m_i], \\
			&z_{t,l}\geq0, &\forall t&\in\mathcal{T},\,l\in\mathcal{L}_{t}.
\end{align}
\end{subequations}

\section{Table of Notation}
\label{Appendix:Nomenclature}

\begin{center}
\footnotesize
	\begin{longtable}{l l}
\toprule
Name & Description \\

\midrule
	\multicolumn{2}{l}{\textbf{GBT Ensemble Definition}} \\
$n$ & Number of the GBT-trained function (continuous) variables \\
$i$ & Continuous variable index \\
$x_i$ & Continuous variable \\
$\bm{x}$ & Vector $(x_1,\ldots,x_n)^T$ \\
$\mathcal{T}$ & Set of gradient boosted trees \\
$t$ & Gradient boosted tree \\
$\mathcal{V}_t$ & Set of split nodes (vertices) in tree $t$ \\
$\mathcal{L}_t$ & Set of leaf nodes in tree $t$ \\
$s$ & Split node associated with a tree $t$ and mainly referred to as $(t,s)$ \\
$i(t,s)$ & Continuous variable index associated with split node $s$ in tree $t$ \\
$v(t,s)$ & Splitting value of variable $x_{i(t,s)}$ at split node $s$ in tree $t$ \\ 
$\text{GBT}_t(\bm{x})$ & Tree $t$ evaluation at point $\bm{x}$ \\
$\text{GBT}(\bm{x})$ & GBT ensemble evaluation at point $\bm{x}$ \\

\midrule
	\multicolumn{2}{l}{\textbf{Convex MINLP with GBTs Problem Definition}} \\
$\text{cvx}(\bm{x})$ & Convex function evaluation at point $\bm{x}$ \\
$m_i$ & Number of variable $x_i$ splitting values \\
$v_{i,j}$ & $j$-th greatest variable $x_i$ splitting value \\
$v_i^L$ or $v_{i,0}$ & Variable $x_i$ lower bound \\
$v_i^U$ or $v_{i,m_i+1}$ & Variable $x_i$ upper bound \\
$\bm{v}^L$ & Vector $(v_1^L,\ldots,v_n^L)$ \\
$\bm{v}^U$ & Vector $(v_1^U,\ldots,v_n^U)$ \\
$\text{Left}_{t,s}$ & Set of leaves in the subtree rooted in the left child of $s$ in tree $t$ \\
$\text{Right}_{t,s}$ & Set of leaves in the subtree rooted in the right child of $s$ in tree $t$ \\
$F_{t,l}$ & Contribution of leaf node $l$ in tree $t$ \\
$y_{i,j}$ & Binary variable indicating whether $x_i\leq v_{i,j}$, or not \\
$z_{t,l}$ & Binary variable specifying whether tree $t$ evaluates at leaf $l$ \\
$d$ & Maximum tree depth \\

\midrule
	\multicolumn{2}{l}{\textbf{Branch-and-Bound Algorithm Overview}} \\
$[\bm{v}^L,\bm{v}^U]$ & Optimization problem global domain \\
$S=[\bm{L},\bm{U}]$ & Optimization problem subdomain / B\&B node \\
$(x_i,v)$ & GBT splitting point / B\&B branch \\
$S_{\text{left}},S_{\text{right}},S_c,S'$ & B\&B nodes \\
$Q$ & Set of unexplored B\&B nodes \\
$P_{\text{root}}$ & Initial GBT ensemble partition at B\&B root node \\
$P,P',P''$ & GBT ensemble partitions \\
$b^{\text{cvx},S}$ & Convex lower bound over domain $S$ \\
	$b^{\text{GBT},S,P}$ & GBT lower bound over domain $S$ with respect to partition $P$ \\

\midrule
	\multicolumn{2}{l}{\textbf{Lower Bounding}} \\
$R^{S}$ & Optimal objective value, i.e., tightest relaxation \\
$\hat{R}^S$ & Relaxation dropping linking constraints \\
$b^{\text{GBT},S,*}$ & Optimal GBT lower bound over domain $S$ \\
$\bm{x}^*$ & Optimal solution \\
$i,j,l$ & Subset indices of a GBT ensemble partition \\
$k$ & GBT ensemble partition size \\
$\mathcal{T}_i,\mathcal{T}_j,\mathcal{T}',\mathcal{T}''$ & Subsets of GBTs \\
$N$ & GBT ensemble subset size \\
	$n^{\mathcal{T}, S}$ & Number of leaves in GBT subset $\mathcal{T}$ over domain $S$\\
	$f^*$ & Best found feasible objective\\
	$q$ & Time limit on lower bound improvement algorithm\\
\midrule
\multicolumn{2}{l}{\textbf{Branching}} \\
$B$ & Branch ordering \\ 
$r((x_i,v),t)$ & Set of nodes in tree $t$ that split on $(x_i,v)$ \\
$d(s)$ & Depth of split node $s$ (root node has zero depth) \\
$w(s)$ & Weight of split node $s$ \\
$i(s)$ & Number of inactive leaves below split $s$ when branching with respect to $(x_i,s)$ \\ 
$\text{weight}((x_i,v),t)$ & Weight assigned to $(x_i,v)$ in tree $t$\\
$\text{weight}((x_i,v),\mathcal{T})$ & Weight assigned to $(x_i,v)$ in GBT ensemble $\mathcal{T}$\\
$\text{inactive}((x_i,v),\mathcal{T})$ & Number of inactive leaves when branching on pair $(x_i,v)$ in $\mathcal{T}$ \\
$\text{cover}(t,s)$ & Set of leaves covered by split node $s$ at tree $t$ \\
$S,S_{\text{left}},S_{\text{right}}, S_0$ & B\&B nodes denoted by their corresponding domain \\
$l$ & Strong branching lookahead parameter\\

\bottomrule
\caption{Nomenclature}
\label{Table:Nomenclature}
\end{longtable}

\end{center}

\section{Heuristics}\label{subsec:heuristics}

\subsection{Mixed-Integer Convex Programming Heuristic}\label{subsubsec:incremental_miqp}
For a given a subset $\mathcal{T}'\subseteq\mathcal{T}$ of trees, let $f_{\mathcal{T}'}(\cdot)$ be the objective function obtained by ignoring the trees $\mathcal{T}\setminus\mathcal{T}'$.
Then, $\min_{\gbtVec{v}^{L}\leq\gbtVec{x}\leq\gbtVec{v}^{U}}\{f_{\mathcal{T}'}(\gbtVec{x})\}$ may be significantly more tractable than the original problem instance when $|\mathcal{T}'|<<|\mathcal{T}|$. 
So, the \cref{alg:miqp_heur} heuristic solves the original convex MINLP by sequentially solving smaller convex MINLP sub-instances of increasing size.
A sub-instance is restricted to a subset $\mathcal{T}'\subseteq\mathcal{T}$ of GBTs.
Let $\mathcal{T}^{(k)}$ be the subset of trees when the $k$-th heuristic iteration begins.
Initially, $\mathcal{T}^{(0)}=\emptyset$, i.e., $f_{\mathcal{T}^{(0)}}(\cdot)$ consists only of the convex part.
Denote by $\gbtVec{x}^{(k)}$ the sub-instance optimal solution minimizing $f_{\mathcal{T}^{(k)}}(\cdot)$.
Note that $\gbtVec{x}^{(k)}$ is feasible for the full instance. 
Each iteration $k$ chooses a set of $N$ additional trees $\mathcal{T}^{\text{next}}\subseteq\mathcal{T}\setminus\mathcal{T}^{(k)}$ and constructs $\mathcal{T}^{(k+1)}=\mathcal{T}^{(k)}\cup\mathcal{T}^{\text{next}}$, i.e., $\mathcal{T}^{(k)}\subseteq\mathcal{T}^{(k+1)}$.
Consider two approaches for picking the $N$ trees between consecutive iterations: (i) training-aware selection and (ii) best improvement selection.
Termination occurs when the time limit is exceeded and \cref{alg:miqp_heur} returns the best computed solution.

\paragraph{Training-aware selection}
Let $T_1,T_2,\ldots,T_m$ be the 
tree generation order during training.
This approach selects the trees $\mathcal{T}^{\text{next}}$ according to this predefined order. 
That is, in the $k$-th iteration, $\mathcal{T}^{(k)}=\{T_1,\ldots,T_{kN}\}$ and $\mathcal{T}^{\text{next}}=\{T_{kN+1},\ldots,T_{(k+1)N}\}$.
A GBT training algorithm constructs the trees iteratively, so each new tree reduces the current GBT ensemble error with respect to the training data.
Thus, we expect that the earliest-generated trees better approximate the learned function than the latest-generated trees.
Specifically, for two subsets $\mathcal{T}_A,\mathcal{T}_B\subseteq\mathcal{T}$ with the property that $t_a<t_b$ for each $T_{t_a}\in\mathcal{T}_A$ and $T_{t_b}\in\mathcal{T}_B$, we expect that $|f_{\mathcal{T}_A}(\gbtVec{x})-f^*(\gbtVec{x})|\leq |f_{\mathcal{T}_B}(\gbtVec{x})-f^*(\gbtVec{x})|$, for each $\gbtVec{v}^{L}\leq\gbtVec{x}\leq\gbtVec{v}^U$, where $f^*$ is the original objective function, i.e., the optimal approximation.
Intuitively, earlier trees place the GBT function within the correct vicinity, while later trees have a fine tuning role.

\paragraph{Best improvement selection}
In this approach, the $k$-th iteration picks the $N$ trees with the maximum contribution when evaluating at $\gbtVec{x}^{(k)}$.
We select $\mathcal{T}^{\text{next}} \subseteq\mathcal{T} \setminus \mathcal{T}^{(k)}$ so that, for each pair of trees $T_t\in\mathcal{T}^{\text{next}}$ and $T_{t'}\in\mathcal{T}\setminus(\mathcal{T}^{(k)}\cup\mathcal{T}^{\text{next}})$, it holds that $f_{t}(\gbtVec{x}^{(k)})\geq f_{t'}(\gbtVec{x}^{(k)})$.
Assuming that approximation $\mathcal{T}^{(k)}$ is poor, then $\mathcal{T}^{\text{next}}$ contains the trees that refute optimality of $\gbtVec{x}^{(k)}$ the most, from the perspective of $f_t(\gbtVec{x}^{(k)})$ $t\in\mathcal{T}\setminus\mathcal{T}^{(k)}$.

\begin{algorithm}
	\caption{Mixed-integer convex programming heuristic}\label{alg:miqp_heur}
	\begin{algorithmic}[1]
	\State $k\gets0$
	\State $\mathcal{T}^{(k)}\gets\emptyset$
	\While{the time limit is not exceeded}
		\State $\gbtVec{x}^{(k)}\gets\argmin\limits_{\gbtVec{v}^L\leq\gbtVec{x}\leq\gbtVec{v}^U} f_{T^{(k)}}(\gbtVec{x})$
		\State{Choose $\mathcal{T}^{\text{next}}$ from $\left\{\mathcal{T}'\,\mid\,\mathcal{T}'\subseteq\mathcal{T}\setminus\mathcal{T}^{(k)}, |\mathcal{T}'|=\min\{N, |\mathcal{T}\setminus\mathcal{T}^{(k)}|\}\right\}$}
		\State$\mathcal{T}^{(k+1)}\gets\mathcal{T}^{(k)}\cup\mathcal{T}^{\text{next}}$
		\State $k\gets k+1$
		\EndWhile
	\State \Return$\argmin\limits_{k\in\{0,\dots, k-1\}}f\left(\gbtVec{x}^{(k)}\right)$
\end{algorithmic}

\end{algorithm}

\subsection{Particle Swarm Optimization}

\cite{488968} introduce PSO for optimizing continuous nonlinear functions.
PSO computes a good heuristic solution by triggering $m$ particles that collaboratively search the feasibility space.
PSO picks the initial particle position $\gbtVec{x}_i^{(0)}$ and search direction $\gbtVec{v}_i^{(0)}$ of particle $i$ randomly.
The search occurs in a sequence of rounds.
In round $k$, every particle chooses its next position $\gbtVec{x}_i^{(k+1)}$ by
following the direction specified by a weighted sum of: (i) the current trajectory direction $\gbtVec{v}_i^{(k)}$, (ii) the particle's best found solution $\gbtVec{p}_i$, (iii) the globally best found solution $\gbtVec{g}$, and moving by a fixed step size.
The \emph{inertia term} $\omega\gbtVec{v_i}^{(k)}$ controls how quickly a particle changes direction.
The \emph{cognitive term} $c_1\cdot r_1\cdot(\gbtVec{p}_i-\gbtVec{x}_i^{(k)})$ controls the particle tendency to move to the best observed solution by that particle.
The \emph{social term} $c_2\cdot r_2\cdot(\gbtVec{g}-\gbtVec{x}_i^{(k)})$ controls the particle tendency to move toward the best solution observed by any particle.
Coefficients $\omega$, $c_1$, and $c_2$ are tunable parameters.
Termination occurs either when all particles are close, or within a specified time limit.	\Cref{alg:pso} lists the PSO algorithm.
\begin{algorithm}[t]
	\caption{Particle Swarm Optimization}\label{alg:pso}
	\begin{algorithmic}
\State Compute initial position $\gbtVec{x_i}^{(0)}\in\mathbb{R}^n$ and velocity $\gbtVec{v_i}^{(0)}\in\mathbb{R}^n$ for each particle $i=1,\ldots,m$.
\State $\gbtVec{p_i}\leftarrow \gbtVec{x_i}^{(0)}$
\State $\gbtVec{g}\leftarrow \text{arg}\min\{f(\gbtVec{p_i})\}$
\State $k\leftarrow 0$
\While {the time limit is not exceeded}
\For {$i=1,\ldots,m$}
\State Choose random values $r_1,r_2\sim U(0,1)$
\State $\gbtVec{v_i}^{(k+1)}\leftarrow\omega\gbtVec{v_i}^{(k)} + c_1\cdot r_1\cdot(\gbtVec{p_i}-\gbtVec{x_i}^{(k)}) + c_2\cdot r_2\cdot(\gbtVec{g}-\gbtVec{x_i}^{(k)})$
\State $\gbtVec{x_i}^{(k+1)}\leftarrow\gbtVec{x_i}^{(k)}+\gbtVec{v_i}^{(k+1)}$
\If {$f(\gbtVec{x_i}^{(k+1)})<f(\gbtVec{p_i})$}
\State $\gbtVec{p_i} \leftarrow \gbtVec{x_i}^{(k+1)}$
\EndIf
\EndFor
\State $\gbtVec{g}\leftarrow \text{arg}\min\{f(\gbtVec{p_i})\}$
\State $k\leftarrow k+1$
\EndWhile
\end{algorithmic}

\end{algorithm}

For Problem~(\ref{eq:full_formulation}), we improve the PSO performance by avoiding initial particle positions in feasible regions strictly dominated by the convex term. 
We project the initial random points close to regions where the GBT term is significant compared to the convex term.

\subsection{Simulated Annealing}\label{app:sa}
\Cref{alg:sa} lists the simulated annealing algorithm \citep{Kirkpatrick671}.
\begin{algorithm}
	\caption{Simulated Annealing}\label{alg:sa}
	\begin{algorithmic}[1]
\State Compute an initial solution $\bm{x}^{(0)}\in\mathbb{R}^n$. 
\State Set initial temperature $T^{(0)}=1$ and probability constant $c=1$.
\State Set temperature factor $\alpha\in[0.80,0.99]$.
\State $t=0$, $k=0$
\While {$T^{(t)}>\epsilon$}
\For {$r$ iterations}
\State Select a neighboring solution $\bm{x}\in\mathcal{N}(\bm{x}^{(k)})$ randomly.
\If {$f(\bm{x})<f(\bm{x}^{(k)})$}
\State $\bm{x}^{(k+1)}\leftarrow\bm{x}$
\State $k\leftarrow k+1$
\Else
\State Choose $p\sim U(0,1)$
\If {$\exp(-(f(\bm{x})-f(\bm{x}^{(k)}))/cT^{(t)})>p$}
\State $\bm{x}^{(k+1)}\leftarrow\bm{x}$
\State $k\leftarrow k+1$
\EndIf
\EndIf
\EndFor
\State $T^{(t+1)}\leftarrow \alpha T^{(t)}$
\State $t\leftarrow t+1$
\EndWhile
\end{algorithmic}

\end{algorithm}
\section{Numerical Results: Heuristic Solutions}\label{app:feasibility}
This section assesses performance of the \cref{subsec:feasibility} heuristic algorithms compared with simulated annealing.
We use CPLEX 12.7 and Gurobi 7.5.2 as: (i) black-box solvers for the entire convex MINLP (\ref{Eq:Optimization_Problem}) and  
(ii) heuristic components for solving convex MINLP (\ref{Eq:Optimization_Problem}) instances in the Section~\ref{subsec:feasibility} convex MINLP heuristic.
The R package GenSA \citep{GenSA} runs the Simulated Annealing (SA) metaheuristic.
We provide a SA technical description \citep{Kirkpatrick671} in \Cref{app:sa}.
The Python module PySwarms \citep{pyswarmsJOSS2018} implements the \Cref{subsec:feasibility} Particle Swarm Optimization (PSO) metaheuristic.
Each heuristic, i.e.\ TA, BI, and Random, uses either CPLEX, or Gurobi as a black-box convex MINLP solver.
We append the labels -C or -G to indicate the underlying solver. 
We use the default CPLEX 12.7 and Gurobi 7.5.2 tolerances, i.e., relative MIP gap, integrality and barrier convergence tolerances of $10^{-4}$, $10^{-5}$ and $10^{-8}$, respectively. 
We use the default SA parameters.
We parameterize PSO with inertia term $\omega=0.5$, cognitive term $c_1=0.7$, social term $c_2=0.3$, 500 particles and an iteration limit of 100.
Each particle takes a randomly generated point, $\gbtVec{x}^{(0)}\in[\gbtVec{v}^L, \gbtVec{v}^U]$, and its projection, $\gbtVec{x}^{(p)}$ on $\gbtMat{P}$ and initializes at $\gbtVec{x}=h\cdot\gbtVec{x}^{(0)} + (1-h)\cdot\gbtVec{x}^{(p)}$.
For our tests, we use $h=0.15$.

\subsection{Concrete Mixture Design}\label{subsubsec:concrete_heur}
\begin{table}
	\TABLE
	{Concrete mixture design instance: black-box solver solutions (upper bounds) by solving the entire mixed-integer convex programming (convex MINLP) model using: (i) CPLEX 12.7, (ii) Gurobi 7.5.2, (iii) Simulated Annealing (SA), and (iv) Particle Swarm Optimization (PSO), with 1 hour timeout.
	\label{tab:concrete_feas_results}}	
	{\begin{tabular}{ccrrrrrrrrr}
	\toprule
	$\lambda$ && CPLEX 12.7 & \hspace*{1cm} & Gurobi 7.5.2 & \hspace*{1cm} & PSO & \hspace*{1cm} & SA \\
	\midrule
	 $1$ && $-14.2$ && $-17.7$ && $-88.7$ && $\bm{-91.3}$ \\
	 $10$ && $422.7$ && $112.6$ && $-86.0$ && $\bm{-86.6}$ \\
	 $100$ && $4,791.6$&& $1,413.7$ && $-80.1$ && $\bm{-80.3}$ \\
	 $1000$ && $48,480.6$ && $14,425.1$ && $\bm{-75.9}$ && $-71.6$ \\
	\bottomrule
\end{tabular}
}
	{}
\end{table}

\begin{figure}
	\centering
	\FIGURE
	{
		\includegraphics[height=23ex,width=0.9\linewidth]{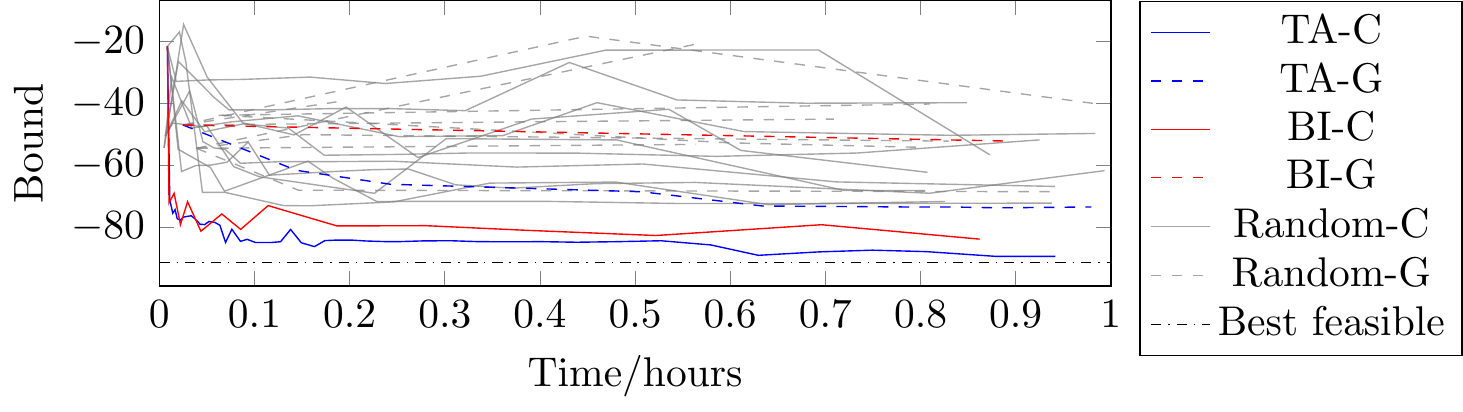}%
	}
	{Concrete mixture design instance ($\lambda=1$): Convex MINLP heuristic using training-aware (TA), best improvement (BI), or random strategies for choosing the next trees. Each iteration selects 10 new trees. The suffixes -C and -G denote using CPLEX 12.7 and Gurobi 7.5.2 as subsolvers, respectively. Best feasible is the simulated annealing solution.
	\label{fig:concrete_heur_compare_lm_1}}
	{}
\end{figure}

\begin{figure}
	\centering
	\FIGURE
	{
		\includegraphics[height=23ex,width=0.9\linewidth]{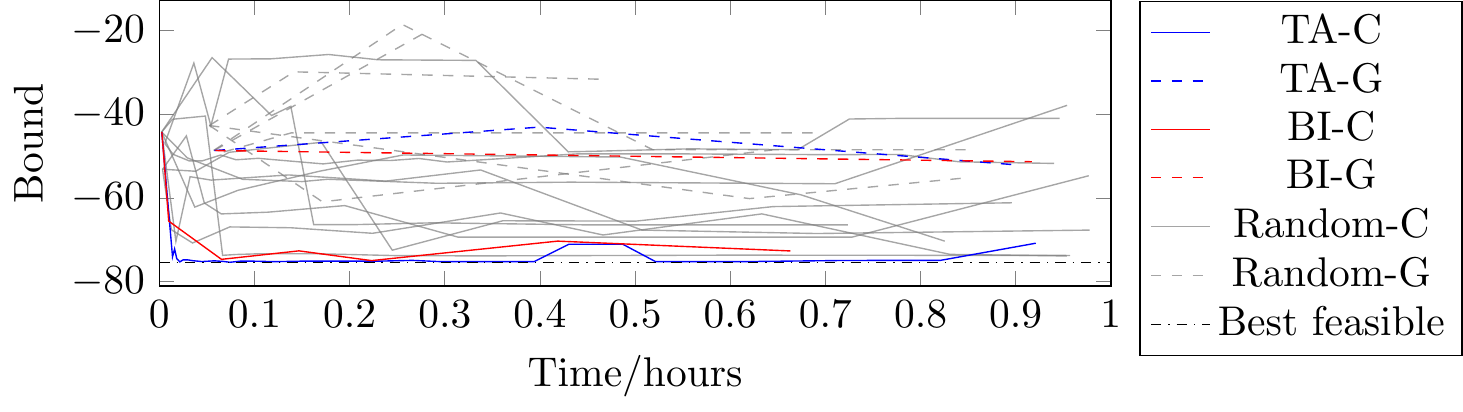}%
	}
	{Concrete mixture design instance ($\lambda=1000$): Convex MINLP heuristic using training-aware (TA), best improvement (BI), or random strategies for choosing the next trees. Each iteration selects 10 new trees. The suffixes -C and -G denote using CPLEX 12.7 and Gurobi 7.5.2 as subsolvers, respectively. Best feasible is the simulated annealing solution.
	\label{fig:concrete_heur_compare_lm_1000}}
	{}
\end{figure}

\Cref{tab:concrete_feas_results} compares the CPLEX 12.7, Gurobi 7.5.2, SA, and PSO computed solutions for the entire convex MINLP, under 1 hour time limit.
SA performs the best.
PSO solution is relatively close to the SA best found solution, compared to CPLEX 12.7 or Gurobi 7.5.2.
\Cref{fig:concrete_heur_compare_lm_1,fig:concrete_heur_compare_lm_1000} evaluate the \cref{subsubsec:incremental_miqp} augmenting convex MINLP heuristic using CPLEX 12.7, Gurobi 7.5.2, and the different tree selection approaches, i.e., (i) training-aware (TA), (ii) best improvement (BI), and (iii) random selection.
\Cref{fig:concrete_heur_compare_lm_1,fig:concrete_heur_compare_lm_1000} also plots the SA best-found solution.
In general, both TA and BI perform better than random selection.
Moreover, TA performs better than BI. 
Therefore, there is a benefit in choosing the earlier trees to find good heuristic solutions.
Interestingly, the solution found in the first iteration of the augmenting convex MINLP heuristic, i.e., by solely minimizing the convex part, is lower than -43, while the upper bounds reported by CPLEX 12.7 and Gurobi 7.5.2 after one hour of execution are greater than -18.

\subsection{Chemical Catalysis}\label{subsubsec:basf_heur}
\Cref{tab:basf_feas_results} compares the CPLEX 12.7, Gurobi 7.5.2, SA, and PSO computed solutions for the entire covex MINLP, under a 1 hour time limit.
SA outperforms all others.
PSO performs well for larger $\lambda$ values,
because it keeps the contribution of the convex part low at initialization.
Gurobi 7.5.2 also performs relatively well for smaller $\lambda$ values, however due to solver tolerances 
	it may report incorrect objective values. For example, using $\lambda=0$ the solver reports an objective of $-174.1$, however a manual evaluation results in $-158.5$.
In fact, both CPLEX 12.7 or Gurobi 7.5.2, may produce incorrect outputs due to solver tolerances, hence a specialized fixing method 
may be necessary.

\begin{table}
	\TABLE
	{Chemical catalysis BASF instance (with different $\lambda$ values): Black-box solver solutions (upper bounds) by solving the entire mixed-integer convex programming (convex MINLP) model using: (i) CPLEX 12.7, (ii) Gurobi 7.5.2, (iii) Simulated Annealing (SA), and (iv) Particle Swarm Optimization (PSO), with 1 hour timeout.\label{tab:basf_feas_results}}	
	{\begin{tabular}{ccrrrrrrrrr}
	\toprule
	$\lambda$ && CPLEX 12.7 && Gurobi 7.5.2 && PSO && SA \\
	\midrule
	$0$    && * & \hspace*{1cm} & $-158.5$ & \hspace*{1cm} & $-96.8$ & \hspace*{1cm} & $\bm{-168.2}$ \\
	$1$    && * && $-101.6$ && $-89.8$ && $\bm{-130.7}$ \\
	$10$   && $952$ && $-100.1$ && $-97.6$ && $\bm{-102.7}$ \\
	$100$  && $1,040$ && $11.5$ && $-82.7$ && $\bm{-84.2}$ \\
	$1000$ && $18,579$ && $606.5$ && $-76.5$ && $\bm{-81.3}$ \\
	\bottomrule
\end{tabular}
}
	{}
\end{table}

\begin{figure}
	\FIGURE
	{
		\includegraphics[height=23ex,width=0.9\linewidth]{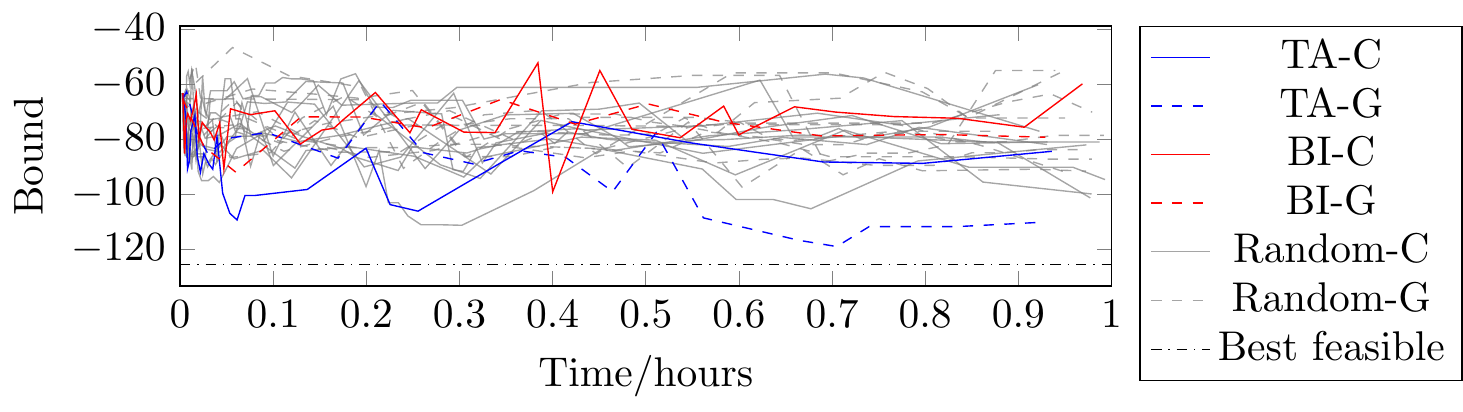}%
	}
	{
		Chemical catalysis instance ($\lambda=1$): Convex MINLP heuristic using training-aware (TA), best improvement (BI), or random strategies for choosing the next trees. Each iteration selects 10 new trees. The suffixes -C and -G denote using CPLEX 12.7 and Gurobi 7.5.2 as subsolvers, respectively. Best feasible is the simulated annealing solution.
	\label{fig:basf_standard_heur_compare_lm_1}}
	{}
\end{figure}

\begin{figure}
	\FIGURE
	{
		\includegraphics[height=23ex,width=0.9\linewidth]{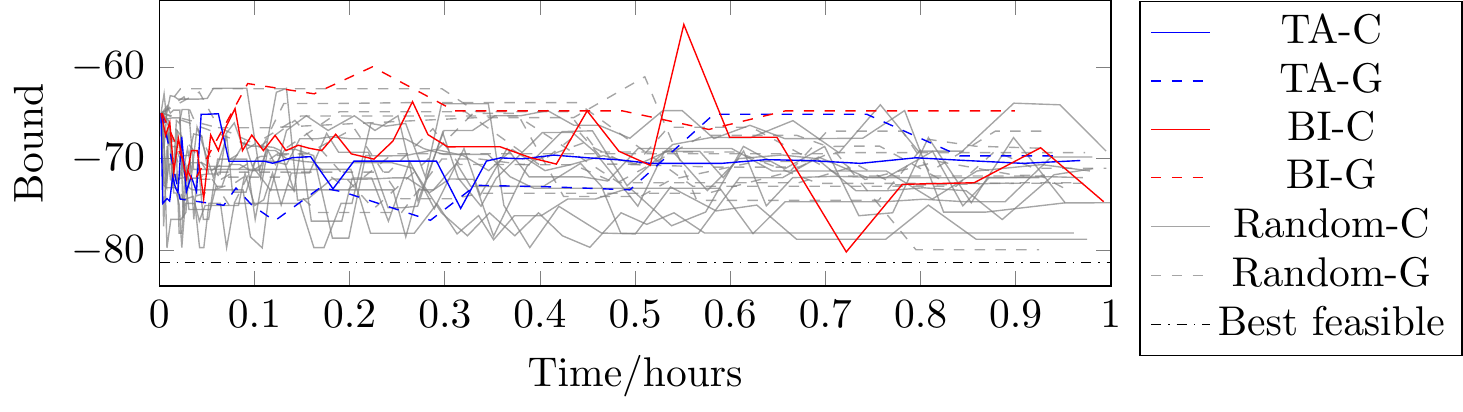}%
	}
	{
		Chemical catalysis instance ($\lambda=1000$): Convex MINLP heuristic using training-aware (TA), best improvement (BI), or random strategies for choosing the next trees. Each iteration selects 10 new trees.  The suffixes -C and -G denote using CPLEX 12.7 and Gurobi 7.5.2 as subsolvers, respectively. Best feasible is the simulated annealing solution.
	\label{fig:basf_standard_heur_compare_lm_1000}}
	{}
\end{figure}

\Cref{fig:basf_standard_heur_compare_lm_1,fig:basf_standard_heur_compare_lm_1000} evaluate the \cref{subsubsec:incremental_miqp} augmenting convex MINLP heuristic for different values of the $\lambda$ input parameter.
We investigate the augmenting convex MINLP heuristic performance using either CPLEX 12.7, or Gurobi 7.5.2 for solving convex MINLP sub-instances and each of the: (i) training-aware (TA), (ii) best improvement (BI), and (iii) random selection strategies.
The \cref{fig:basf_standard_heur_compare_lm_1,fig:basf_standard_heur_compare_lm_1000} best feasible solution is the one produced by SA.
For $\lambda=1$, TA constructs several heuristic solutions that outperform both the BI and random selection ones. 
In this case,  since the GBT part dominates the convex part, TA iteratively computes a better GBT approximation.
For $\lambda=1000$, TA and BI exhibit comparable performance, with BI finding the best solution.
Random selection also performs well because the convex part dominates the GBT part.

\end{APPENDICES}

\end{document}